\documentclass[fleqn,reqno,11pt,a4paper,final]{amsart}
\usepackage[margin=1.9cm]{geometry}
\usepackage{amsmath}
\usepackage{amssymb}
\usepackage{amsthm}
\usepackage{amscd}
\usepackage[ansinew]{inputenc}
\usepackage{cite}
\usepackage{bbm}
\usepackage[usenames, dvipsnames]{color}
\usepackage[english=american]{csquotes}
\usepackage[final]{graphicx}
\usepackage{hyperref}
\usepackage{calc}
\usepackage{verbatim}
\usepackage{enumerate}

\numberwithin{equation}{section}

\newtheoremstyle{thmlemcorr}{10pt}{10pt}{\itshape}{}{\bfseries}{.}{10pt}{{\thmname{#1}\thmnumber{ #2}\thmnote{ (#3)}}}
\newtheoremstyle{thmlemcorr*}{10pt}{10pt}{\itshape}{}{\bfseries}{.}\newline{{\thmname{#1}\thmnumber{ #2}\thmnote{ (#3)}}}
\newtheoremstyle{remexample}{10pt}{10pt}{}{}{\bfseries}{.}{10pt}{{\thmname{#1}\thmnumber{ #2}\thmnote{ (#3)}}}

\theoremstyle{thmlemcorr}
\newtheorem{theorem}{Theorem}
\numberwithin{theorem}{section}
\newtheorem{lemma}[theorem]{Lemma}

\theoremstyle{thmlemcorr*}
\newtheorem{theorem*}{Theorem}
\newtheorem{lemma*}[theorem]{Lemma}
\newtheorem{corollary*}[theorem]{Corollary}
\newtheorem{proposition*}[theorem]{Proposition}
\newtheorem{problem*}[theorem]{Problem}
\newtheorem{conjecture*}[theorem]{Conjecture}
\newtheorem{definition*}[theorem]{Definition}

\theoremstyle{remexample}
\newtheorem{remark}[theorem]{Remark}
\newtheorem{example}[theorem]{Example}

\graphicspath{{../Pictures/}}

\newcommand{\Crm}{\mathrm{C}}

\newcommand{\Lrm}{\mathrm{L}}

\newcommand{\Wrm}{\mathrm{W}}
\newcommand{\Xrm}{\mathrm{X}}

\newcommand{\Fcal}{\mathcal{F}}

\newcommand{\Lcal}{\mathcal{L}}

\newcommand{\Rcal}{\mathcal{R}}

\newcommand{\Wcal}{\mathcal{W}}

\newcommand{\Abb}{\mathbb{A}}

\newcommand{\Nbb}{\mathbb{N}}

\newcommand{\Pbb}{\mathbb{P}}

\DeclareMathOperator{\diam}{diam}

\DeclareMathOperator{\diverg}{div}

\newcommand{\ee}{\mathrm{e}}

\newcommand{\setb}[2]{\bigl\{\, #1 \ \ \textup{\textbf{:}}\ \ #2 \,\bigr\}}

\newcommand{\norm}[1]{\|#1\|}

\newcommand{\abs}[1]{|#1|}

\newcommand{\absBB}[1]{\biggl|#1\biggr|}

\newcommand{\dpr}[1]{\langle #1 \rangle}

\newcommand{\dprb}[1]{\bigl\langle #1 \bigr\rangle}

\newcommand{\cl}[1]{\overline{#1}}
\newcommand{\di}{\mathrm{d}}
\newcommand{\dd}{\;\mathrm{d}}
\newcommand{\DD}{\mathrm{D}}
\newcommand{\N}{\mathbb{N}}
\newcommand{\R}{\mathbb{R}}

\newcommand{\ONE}{\mathbbm{1}}

\newcommand{\toweak}{\rightharpoonup}
\newcommand{\toweakstar}{\overset{*}\rightharpoonup}

\newcommand{\todown}{\downarrow}

\newcommand{\embed}{\hookrightarrow}
\newcommand{\cembed}{\overset{c}{\embed}}

\newcommand{\sbullet}{\begin{picture}(1,1)(-0.5,-2.5)\circle*{2}\end{picture}}
\newcommand{\frarg}{\,\sbullet\,}

\newcommand{\eps}{\epsilon}

\newcommand{\meantmp}[2]{#1\langle{#2}#1\rangle}
\newcommand{\mean}[1]{\meantmp{}{#1}}

\providecommand{\skptmp}[3]{{\ensuremath{#1\langle {#2}, {#3} #1\rangle}}}
\providecommand{\skp}[2]{\skptmp{}{#1}{#2}}

\newcounter{assumption}
\makeatletter
\newcommand{\nextas}[1]{%
 ~\refstepcounter{assumption}%
   \protected@write \@auxout{}{\string\newlabel{#1}{{(A\theassumption)}{\thepage}{(A\theassumption)}{#1}{}}}%
   \hypertarget{#1}{(A\theassumption)}%
}
\makeatother

\def\Xint#1{\mathchoice
{\XXint\displaystyle\textstyle{#1}}%
{\XXint\textstyle\scriptstyle{#1}}%
{\XXint\scriptstyle\scriptscriptstyle{#1}}%
{\XXint\scriptscriptstyle\scriptscriptstyle{#1}}%
\!\int}
\def\XXint#1#2#3{{\setbox0=\hbox{$#1{#2#3}{\int}$}
\vcenter{\hbox{$#2#3$}}\kern-.5\wd0}}
\def\dashint{\,\Xint-}

\renewcommand{\eps}{\varepsilon}
\renewcommand{\epsilon}{\varepsilon}
\renewcommand{\phi}{\varphi}
\renewcommand{\hat}{\widehat}

\begin{document}


\title[Regularity and approximation of rate-independent systems]{Regularity and approximation of strong solutions to rate-independent systems}

\author{Filip Rindler}
\address{\textit{Filip Rindler:} Mathematics Institute, University of Warwick, Coventry CV4 7AL, UK}
\email{F.Rindler@warwick.ac.uk}

\author{Sebastian Schwarzacher}
\address{\textit{Sebastian Schwarzacher:} Department of Analysis, Charles University Prague, Sokolovsk\'{a} 83, 186 75 Praha 8, Czech Republic}
\email{schwarz@karlin.mff.cuni.cz}

\author{Endre S\"uli}
\address{\textit{Endre S\"{u}li:} Mathematical Institute, University of Oxford, Andrew Wiles Building, Woodstock Road, Oxford OX2 6GG, United Kingdom}
\email{endre.suli@maths.ox.ac.uk}

\maketitle
\thispagestyle{empty}

\begin{abstract}
Rate-independent systems arise in a number of applications. Usually, weak solutions to such problems with potentially very low regularity are considered, requiring mathematical techniques capable of handling nonsmooth functions. In this work we prove the existence of H\"{o}lder-regular strong solutions for a class of rate-independent systems. We also establish additional higher regularity results that guarantee the uniqueness of strong solutions. The proof proceeds via a time-discrete Rothe approximation and careful elliptic regularity estimates depending in a quantitative way on the (local) convexity of the potential featuring in the model. In the second part of the paper we show that our strong solutions may be approximated by a fully discrete numerical scheme based on a spatial finite element discretization, whose rate of convergence is consistent with the regularity of strong solutions whose existence and uniqueness are established.

\vspace{4pt}



\vspace{4pt}

\noindent\textsc{Date:} \today{}. 
\end{abstract}

\setcounter{tocdepth}{2}

\section{Introduction}

Rate-independent systems are used to model a plethora of physical phenomena for which the speed of evolution does not influence the extent of energy dissipation, including elasto-plasticity, damage and delamination in solids, crack propagation, and shape-memory alloys. We refer to the recent monograph~\cite{MielkeRoubicek15book} for an up-to-date overview of the literature on both the theory and the applications of rate-independent systems. In this work we consider only \emph{purely dissipative} systems, i.e.\ those without elastic variables (which, by definition, are those quantities that can be changed without dissipating energy). A formal prototype of a rate-independent system of this kind is
\[
 \frac{\dot{u}}{\abs{\dot{u}}} - \Delta u + \DD W_0(u) = f, \qquad u \colon [0,T] \times \Omega \to \R^m,
\]
where $W_0$ is a smooth potential defined on $\R^m$ and $f \colon [0,T] \times \Omega \to \R^m$ is a given external force.
For the mathematical analysis, the crucial feature of such systems of equations is the quasi-static nature of their evolution, namely that solutions simply \emph{rescale} when scaling the external forces. In a sense, quasi-static evolutions \emph{follow} the total (energetic and external) forces with infinite speed and hence the evolution should, at least in sufficiently convex situations, \enquote{inherit} regularity from the external force. Questions about the regularity of solutions were in fact already raised in the early stages of development of the modern energetic theory of rate-independent systems; see in particular~\cite{MielkeTheilLevitas02} and~\cite{MielkeTheil04} (Section~7.3 in the second reference discusses temporal regularity for the uniformly convex case). Further investigations in that direction were carried out recently in~\cite{ThomasMielke10,MielkeZelik14}.

The theory of rate-independent systems is dominated by several notions of weak solution concepts (see~\cite{Mielke11,MielkeRoubicek15book}), and the question how \enquote{realistic} particular notions of solution are from the point of view of applications has been the source of considerable discussion in the literature. A central question is whether there is \enquote{foresight} in the system, i.e.\ whether the system will dissipatively jump over potential energy barriers in order to realize an ultimate drop in energy. The Mielke--Theil energetic solutions do this, but more recent concepts, such as bounded-viscosity solutions~\cite{MiRoSa09MSJR,MiRoSa12BVSV}, do not. A more substantial discussion of this point can be found in~\cite{MielkeRoubicek15book} and Section~\ref{ssc:ex} of the present paper.

Our objective here is to show that in a \enquote{mildly nonconvex} setting
with a quadratic gradient regularizer, one may focus right from the start on \emph{strong solutions},
thereby avoiding much of the ambiguity. We will show that:
\begin{itemize}
\item strong solutions exist;
\item they have additional regularity and satisfy quantitative H\"{o}lder estimates in space and time;
\item they are unique within their natural regularity class;
\item and they can be numerically approximated by a convergent finite element method.
\end{itemize}
Thus, in this setting, there is no need to consider weak solution concepts. We refer to Section~\ref{ssc:ex} for a motivating scalar example, from which it can be seen that weak solutions are nonunique, jump \enquote{too early}, and do not agree with strong solutions, even if both exist.  This motivates our focus on strong solutions (where they exist), which seem to be the physically realistic ones.

Our assumptions apply for instance to some double-well energies with \enquote{mild} nonconvexities, where \enquote{mild} refers to the fact that the admissible extent of nonconvexity depends on the $\Lrm^2$-Poincar\'{e} embedding constant of the domain. The class of problems considered here involves an elliptic regularization, which ensures that the \emph{regularized} functional is convex. This framework is frequently used in the literature, as some regularization is already needed in order to prove the existence of solutions; see~\cite{MielkeRoubicek15book}.
Indeed, in the model considered here, the regularity of solutions is related to the (local) convexity of the energy functional associated with the model. This highlights that the assumed (non-)convexity has to be appropriate in order for one to be able to show the existence, uniqueness and regularity of solutions. It is for this reason that the existence, uniqueness and regularity theory developed in this paper remains valid up until the moment when the solution reaches a point of
non-convexity; see Remark~\ref{rem:H}.

Partial results in the direction of regularity have already been presented in~\cite{ThomasMielke10,MielkePaoliPetrovStefanelli10,MielkeZelik14,KneesRossiZanini15}, but those papers deal with the Mielke--Theil energetic solution concept. Here, we are able to guarantee more regularity in both time and space, and with the natural dependence of the regularity of the solution on that of the data. We also establish quantitative estimates and develop a direct, unified approach, avoiding abstract hypotheses.
Some remarks on uniqueness in this framework can be found in Sections 3.4.3 and 3.4.4 in~\cite{MielkeRoubicek15book}.

As for the numerical approximation of rate-independent systems,
one concrete application was considered by Han and Reddy \cite{HR00}, who provide the mathematical analysis
of semidiscrete and fully discrete approximations to the primal problem in elastoplasticity.
They proved the strong convergence of the sequence of numerical approximations defined by their method(s)
under conditions of minimal regularity of the solution, that is, without any assumptions on the regularity of
the solution over and above those established in the proof of well-posedness of the problem. More recently,
Mielke et al. \cite{MielkePaoliPetrovStefanelli10} have developed error estimates for space-time
discretizations in the context of evolutionary variational inequalities of rate-independent type.
After introducing a general abstract evolution problem, they formulated a fully discrete approximation
and developed a-priori error estimates. Bartels \cite{B14} derived quasioptimal error bounds for implicit
discretizations of a class of rate-independent evolution problems, without imposing regularity assumptions
on the exact solutions, but the load functionals involved in the problem were assumed to be twice
continuously differentiable in time. The new contribution in this paper to the mathematical analysis of finite element
approximations of rate-independent systems is that the orders of convergence established in Section \ref{sec:num} are
in agreement with the regularity results for the analytical solution developed in Section \ref{sc:estimate}.

\subsection{Strong solutions}

For technical reasons we only consider the physically relevant case of $\Omega \subset \R^d$ being a bounded open set, with $d \in \{2,3\}$. Some of our results are also true in higher dimensions, but additional (more restrictive) assumptions are then necessary.

Specifically, we will investigate the class of rate-independent systems that can be formulated as follows: for $\Omega \subset \R^d$ a bounded open Lipschitz domain, $d=2,3$, and $T > 0$, consider the following (formal) system:
\begin{align} \label{eq:PDE}
  \left\{
  \begin{aligned}
    \partial \Rcal_1(\dot{u}(t)) - \Lcal_t u(t) + \DD\Wcal_0(u(t)) &\ni f
      \quad \text{in $(0,T] \times \Omega$,}\\
    u(t)|_{\partial \Omega} &= 0  \quad\text{for $t \in (0,T]$,} \\
    \qquad u(0) &= u_0.
  \end{aligned} \right.
\end{align}
Here, $\Rcal_1 \colon \Lrm^1(\Omega;\R^m) \to \R \cup \{+\infty\}$ is the \textit{rate-independent dissipation potential}, which is assumed to be proper (in the sense that it is not identically $+\infty$), convex, and positively $1$-homogeneous; $\partial \Rcal_1$ is its subdifferential;
$f \in \Wrm^{1,a}(0,T;\Lrm^p(\Omega;\R^m))$, for $a \in (1,\infty)$ and $p \in [2,\infty)$, is
the \textit{external loading (force)}, $u_0 \in (\Wrm^{1,2}_0 \cap \Lrm^q)(\Omega;\R^m)$,
$q \in (1,\infty)$ is the \textit{initial value}, and $\Wcal \colon \Lrm^q(\Omega;\R^m) \to \R$
denotes the \textit{(elastic) energy functional}; we
assume that $\DD \Wcal_0(v) \in \Lrm^{q/(q-1)}(\Omega;\R^m)$ and
\[\norm{\DD \Wcal_0(v)}_{\Lrm^{q/(q-1)}}
\leq C(1+\norm{v}_{\Lrm^q}^{q-1})\qquad \mbox{for all $v \in \Lrm^q(\Omega;\R^m)$ and some $q \in (1,\infty)$}.
\]
Finally, the \textit{regularizer} $\Lcal_t$ is a (possibly time-dependent) second-order linear partial
differential operator in the spatial variables (most commonly, $\Lcal_t = \Delta$).
Our precise assumptions on $\Lcal_t$ are detailed below.

We note that our choice of a homogeneous Dirichlet boundary condition is made merely to simplify the exposition; analogous results for nonhomogeneous Dirichlet boundary data necessitate only unessential technical alterations to our arguments.

We call a map $u \in \Lrm^\infty(0,T;(\Wrm^{1,2}_0 \cap \Lrm^\infty)(\Omega;\R^m))$, such that its weak time derivative $\dot{u}$ has regularity $\dot{u} \in \Lrm^1(0,T;\Wrm^{1,2}_0(\Omega;\R^m))$, a \textit{strong solution} to~\eqref{eq:PDE} provided that
\begin{equation} \label{eq:strongsol_dualreg}
  \Lcal_t u(t) - \DD\Wcal_0(u(t)) + f(t) \in \Lrm^\infty(0,T;\Lrm^2(\Omega;\R^m))
\end{equation}
and
\[
  \left\{ \begin{aligned}
    \Lcal_t u(t) - \DD\Wcal_0(u(t)) + f(t) &\in \partial \Rcal_1(\dot{u}(t))  \qquad\text{for a.e.\ $t \in (0,T]$,} \\
    u(0) &= u_0.
  \end{aligned} \right.
\]
By recalling the definition of the subdifferential and using
$\langle \cdot , \cdot \rangle$ to denote the duality pairing between
$(\Wrm^{1,2}_0(\Omega;\mathbb{R}^m))'$ and $\Wrm^{1,2}_0(\Omega;\mathbb{R}^m)$ with $\Lrm^2(\Omega;\R^m)$ as pivot space, where $(\Wrm^{1,2}_0(\Omega;\mathbb{R}^m))'$ denotes
the dual space of $\Wrm^{1,2}_0(\Omega;\mathbb{R}^m)$,
the above means that
\begin{equation}\label{eq:strongsol_ineq}
\left\{ \begin{aligned}
  &\Rcal_1(\dot{u}(t)) + \dprb{\Lcal_t u(t) - \DD\Wcal_0(u(t)) + f(t), \xi(t) - \dot{u}(t)} \leq \Rcal_1(\xi(t))\\
  &\hspace{60pt}\text{for all $\xi \in \Lrm^1(0,T;\Wrm^{1,2}_0(\Omega;\R^m))$ and a.e.\ $t \in (0,T]$},\\
  &\hspace{8cm}u(0) = u_0,
\end{aligned}
\right.
\end{equation}
with the implicit understanding that $\dot{u}(t)$ belongs to the \textit{effective domain}
\[
  D(\mathcal{R}_1):=\setb{ v \in \Lrm^1(\Omega; \R^m)}{ \mathcal{R}_1(v)<+\infty }
\]
of the convex mapping $\mathcal{R}_1$; our assumptions on $\mathcal{R}_1$ stated in~(A\ref{as:second}) below, together with the hypotheses on the initial data formulated in (A\ref{as:last}), will ensure that this is indeed the case. We note here that the attainment of the initial condition makes sense since our regularity assumptions on $u, \dot{u}$ imply that $u \in \Crm([0,T];\Lrm^2(\Omega;\R^m))$.

\begin{remark} \label{rem:exponents}
The above regularity class for $u,\dot{u}$ is not minimal from the point of view of ensuring that~\eqref{eq:strongsol_ineq} is meaningful, but it turns out to be a suitable class of functions for establishing both the existence and the uniqueness of a strong solution.
An alternative concept of solution would be to require that
\[\mbox{$u \in \Lrm^\infty(0,T;(\Wrm^{1,2}_0 \cap \Lrm^q)(\Omega;\R^m))\qquad$ with $\qquad\dot{u} \in \Lrm^1(0,T;(\Wrm^{1,2}_0 \cap \Lrm^q)(\Omega;\R^m))$}\]
such that~\eqref{eq:strongsol_ineq} is satisfied for all $\xi \in \Lrm^1(0,T;(\Wrm^{1,2}_0 \cap \Lrm^q)(\Omega;\R^m))$ (with $\dpr{\frarg,\frarg}$ now also denoting the $(\Lrm^q \times \Lrm^{q/(q-1)})$-duality pairing). Note that here we do not assume~\eqref{eq:strongsol_dualreg} as part of the solution concept.
If $d = 3$ and $1 \leq q \leq 6$, by the Sobolev embedding theorem we have that $(\Wrm^{1,2}_0 \cap \Lrm^q)(\Omega;\R^m) = \Wrm^{1,2}_0(\Omega;\R^m)$, so that the solution we construct in Theorem~\ref{thm:exist} below is also a solution in this sense. However, if $q > 6$, then our existence proof does not provide a solution in the modified sense. Still, uniqueness holds in this larger class, as will follow from the proof of Theorem~\ref{thm:unique}, also see Remark~\ref{eq:strongsol_weak}.
\end{remark}

\subsection{Motivating example} \label{ssc:ex}

Our aim in this section is to motivate and discuss why the existence and regularity of strong solutions should not merely be seen as a theory of \enquote{improved regularity} for weak solutions (for instance in the energetic sense), and why strong solutions should be considered from the outset as a suitable solution concept. To see this, we consider the following zero-dimensional double-well setup:
\[
  \Wcal_0(z) = W_0(z) := \min\{ z(z+2), z(z-2) \},  \qquad
  \Rcal_1(z) = R_1(z) := \abs{z}.
\]
For
\[
  f(t) := t,  \qquad t \in [0,\infty),
\]
an energetic weak solution $u^{\mathrm{weak}}$ with initial value $u(0) = -1$ exists on the infinite time interval $[0,\infty)$. Indeed, it can be checked in an elementary fashion that
\[
  u^{\mathrm{weak}}(t) := \begin{cases}
    -1   & \text{if $t \in [0,1)$,} \\
    \frac{t+1}{2}  & \text{if $t \in [1,\infty)$}
  \end{cases}
\]
satisfies the \emph{energy balance}
\begin{equation} \label{eq:energy_balance}
  E(t,u(t)) - E(0,u(0)) = \int_0^t f'(s) u(s) \dd s - \mathrm{Var}_{R_1}(u;[0,t)),
\end{equation}
where $E(t,z) := W(z) - f(t)z$ for all $t \in [0,\infty)$, and $\mathrm{Var}_{R_1}$ denotes the total $R_1$-variation, as well as the \emph{global stability} inequality
\[
  E(t,u(t)) \leq E(t,z) + R_1(z-u(t))  \quad
  \text{for all $z \in \R$.}
\]

On the other hand, over the time interval $[0,3)$ also a \emph{strong solution} $u^{\mathrm{strong}}$ exists, namely
\[
  u^{\mathrm{strong}}(t) := \begin{cases}
    -1   & \text{if $t \in [0,1)$,} \\
    \frac{t-3}{2}  & \text{if $t \in [1,3)$.}
  \end{cases}
\]
We see that $u^{\mathrm{weak}}$ and $u^{\mathrm{strong}}$ coincide on the interval $[0,1)$, but they differ on the interval $[1,3)$ where both a weak and a strong solution exist. This fact is closely related to the nonuniqueness of weak solutions: both $u^{\mathrm{weak}}$ and $u^{\mathrm{strong}}$ are in fact weak solutions (as can be easily verified), but only one of them is strong. However, $u^{\mathrm{weak}}$ jumps \enquote{too early} and should be regarded as unphysical since it implies \enquote{foresight} in the system, i.e.\ the system jumps as early as possible to minimize $z \mapsto E(t,z) + R_1(z-u(t))$ at all times. These facts are of course well-known in the theory of rate-independent systems (see e.g.~\cite{MielkeRoubicek15book} and also~\cite{MielkePaoliPetrovStefanelli10}).

Furthermore, $u^{\mathrm{strong}}$ can be extended to $t \in [3,\infty)$, as a weak solution only, by setting
\[
  u^{\mathrm{ext}}(t) := \begin{cases}
    -1   & \text{if $t \in [0,1)$,} \\
    \frac{t-3}{2}  & \text{if $t \in [1,3)$,} \\
    \frac{t+1}{2}  & \text{if $t \in [3,\infty)$.}
  \end{cases}
\]
This function satisfies the energy balance~\eqref{eq:energy_balance} as well as the \emph{local stability} condition
\[
  -\DD E(t,u(t)) \in \partial R_1(0).
\]
We remark that $u^{\mathrm{ext}}$ can be understood as a strong solution in $[0,3)$ and $(3,\infty)$ with regard to the \emph{convex} potentials $z \mapsto z(z+2)$ and $z \mapsto z(z-2)$, respectively. Only at the jump point $t = 3$ do we need to use the concept of energetic solution. In fact, our existence and uniqueness results, Theorems~\ref{thm:exist} and~\ref{thm:unique}, can be used to show that on $[0,3) \cup (3,\infty)$ the strong solution exists (which we already know from our explicit construction) and that it is unique within the class of strong solutions.

We conclude from this discussion that a physically realistic notion of solution for a rate-independent system could be that of a \enquote{maximally-strong} energetic solution, that is, a Mielke--Theil energetic solution such that on every open interval $I \subset [0,T]$ where it agrees with a strong solution, this strong solution cannot be extended beyond $I$. Clearly, $u^{\mathrm{ext}}$ satisfies this definition and the jump in $u^{\mathrm{ext}}$ occurs precisely when the solution moves between regions of convexity of $E$.

In higher dimensions, this definition has to be refined, which we postpone to future work. In the present paper we aim to prepare this approach by establishing maximal regularity, uniqueness, and approximation properties of strong solutions.

Related examples can be found in Section~1.8 of~\cite{MielkeRoubicek15book}. See also~\cite{Rou15}, where the relationship between conventional weak solutions and local solutions is discussed, and a suitably integrated maximal-dissipation principle is devised to select force-driven local solutions and eliminate solutions with ``too-early jumps'', which may otherwise arise if the notion of solution is simply energy-driven.

\subsection{Assumptions} \label{sc:assume}

Unless otherwise stated, the following conditions will be assumed to hold in the rest of the paper:

\begin{enumerate}[({A}1)]
  \item \label{as:first} Let $d\in\{2,3\}$ and let $\Omega \subset \R^d$ be open, bounded and with boundary of class $\Crm^{1,1}$.
  \item \label{as:second} The rate-independent dissipation (pseudo)potential $\Rcal_1 \colon \Lrm^1(\Omega;\R^m) \to \R \cup \{+\infty\}$ is given as
\[\qquad
  \Rcal_1(v) = \int_\Omega R_1(v(x)) \dd x,  \qquad v \in \Lrm^1(\Omega;\R^m),
\]
with $R_1 \colon \R^m \to \R$ convex, lower semicontinuous, and positively $1$-homogeneous, i.e.\ $R_1(\alpha w) = \alpha R_1(w)$ for any $\alpha \geq 0$ and $w \in \R^m$. As any convex function defined on $\R^m$ is locally Lipschitz continuous, the assumed positive $1$-homogeneity of $R_1$ then implies that $R_1$ is in fact globally Lipschitz continuous on $\R^m$.

  \item \label{as:funct} The energy functional $\Wcal_0 \colon \Lrm^q(\Omega;\R^m) \to [0,\infty]$, where $q \in (1,\infty)$, has the form
\[\qquad
  \Wcal_0(u) = \int_\Omega W_0(u(x)) \dd x,
\]
with $W_0 \in \Crm^1(\R^m;[0,\infty))$ satisfying the following assumptions for suitable constants $C,\mu > 0$, and all
$v, w \in \R^m$:
\begin{align}
  \qquad
  C^{-1}(\abs{v}^q-1) &\leq W_0(v) \leq C(\abs{v}^q+1);  \label{eq:Wgrowth}\\
  \abs{\DD W_0(v)}  &\leq C(1+\abs{v}^{q-1});               \label{eq:DWgrowth}\\
 -\mu \abs{v-w}^2 &\leq (\DD W_0(v)-\DD W_0(w))\cdot(v-w);  \label{eq:Wmon}\\
  \mu C_P(\Omega)^2 &< \kappa.                                    \label{eq:muCP}
\end{align}
Here, $\kappa$ is defined via the ellipticity of the regularizer below, and $C_P(\Omega) > 0$ is the (best) $\Lrm^2$-Poincar\'{e} embedding constant of $\Omega$, i.e.\ the smallest constant $C > 0$ such that
\[\qquad
  \norm{v}_{\Lrm^2} \leq C \norm{\nabla v}_{\Lrm^2}  \qquad
  \text{for all $v \in \Wrm^{1,2}_0(\Omega;\R^m)$.}
\]

 \item The regularizer $\Lcal_t$ is a second-order linear partial differential operator of the form
\[ \qquad
  [\Lcal_t v]^{\beta}(x) := \diverg_x [ \Abb(t,x) : \nabla v(x) ]^{\beta} = \sum_j \partial_j \sum_{\alpha,i} A^{\alpha, \beta}_{i,j}(t,x) \, \partial_i v^\alpha(x)
\]
for $1 \leq \alpha, \beta \leq m$, $1 \leq i,j \leq d$. We assume that the fourth-order tensor $\Abb = (A^{\alpha,\beta}_{i,j})^{\alpha,\beta}_{i,j}$ satisfies the following continuity, ellipticity and symmetry conditions:
\begin{align} \qquad
  A^{\alpha,\beta}_{i,j} &= A^{\beta,\alpha}_{j,i} \in \Crm^{0,1}([0,T]\times\overline{\Omega})
    \quad\text{for $\alpha,\beta\in \{1,\ldots,N\}$, $i,j\in\{1,\ldots,d\}$;} \label{eq:Lcoeff} \\
  \hspace{-10pt} \xi : \Abb(t,x) : \xi &= \sum_{i,j,\alpha,\beta} \xi^\alpha_i A^{\alpha, \beta}_{i,j}(t,x) \xi^\beta_j \geq \kappa \abs{\xi}^2
    \quad\text{for $\xi \in \R^{m \times d}$, $(t,x) \in [0,T]\times \overline\Omega$,}  \label{eq:Lell}
\end{align}
where $\kappa > 0$ is the \emph{modulus of ellipticity}.

  \item \label{as:f} The external force has regularity $f \in \Wrm^{1,a}(0,T;\Lrm^p(\Omega;\R^m))$, for some $a \in (1,\infty]$ and $p \in [2,\infty]$.

  \item \label{as:last} The initial value satisfies $u_0 \in (\Wrm^{1,2}_0 \cap \Lrm^q)(\Omega;\R^m)$ and $u_0$ is \emph{stable}, i.e.\ $u_0$ minimizes the functional
\[
  v \mapsto \int_\Omega R_1(v - u_0) + \nabla v : \frac{\Abb(0,\frarg)}{2} : \nabla v + W_0(v) - f(0) \cdot v \dd x,
\]
where $v \in (\Wrm^{1,2}_0 \cap \Lrm^q)(\Omega;\R^m)$.

  \item \label{as:unique} For the uniqueness result stated in Theorem~\ref{thm:unique}
  and the convergence result in Theorem~\ref{thm:approx}, we also assume that there is a $\sigma\in (0,1)$ such that for every $K > 0$ there exists an $M > 0$, such that
\begin{align}
\label{eq:DDDW} \qquad
\abs{\DD^3 W_0(v)}&\leq\max\{\abs{v}^{-\sigma},M\}\quad \text{ for }\abs{v}\leq K.
\end{align}
\end{enumerate}

These assumptions are for instance satisfied by the nonconvex double-well potential
\[
W_0(v) = \gamma (|v|^2-1)^2\qquad \mbox{for $0<\gamma<\kappa/(4C_P(\Omega))^2$ \;(and $q = 4$, $\mu = 4\gamma$).}
\]
Further examples of potentials satisfying the above hypotheses are
\[ W_0(v)=|v|^q \qquad \mbox{for $q \geq 2$},\qquad\mbox{and}\qquad
W_0(v)=(|v|+1)^{q-2}\abs{v}^2 \qquad \mbox{for $q >1 $}.\]

\begin{remark}
We note that we are allowing nonconvexity in $W_0$, but, by~\eqref{eq:Wmon} and~\eqref{eq:muCP}, this nonconvexity cannot be too strong. It can be shown that our conditions imply that the combined energy functional $\Wcal(u) :=  \int_\Omega \nabla u:\frac{\Abb(t,\cdot)}{2}:\nabla u\dd x + \Wcal_0(u)$ is convex. As a matter of fact, a certain amount of convexity of $W_0$ is necessary, since for strongly nonconvex $W_0$ counterexamples to regularity exist, see the example above or~\cite{MielkeTheil04}.
\end{remark}

\subsection{Main results}

In the course of this work we will prove the following result concerning the existence,
uniqueness, and regularity of a strong solution to~\eqref{eq:PDE}.

\begin{theorem} \label{thm:exist}
Under the assumptions~(A\ref{as:first})--(A\ref{as:last}), there exists a strong solution
\[ \mbox{$u \in \Lrm^\infty(0,T;(\Wrm^{1,2}_0 \cap \Lrm^\infty)(\Omega;\R^m))\qquad$ with
 $\qquad\dot{u} \in \Lrm^1(0,T;\Wrm^{1,2}_0(\Omega;\R^m))$},\]
such that~\eqref{eq:strongsol_dualreg} holds, and $u$ has the following
additional regularity properties:
\begin{enumerate}[(i)]
\item $\nabla^2 u\in \Lrm^\infty (0,T;\Lrm^p(\Omega;\R^{m \times d \times d}))$;
\item $\nabla \dot{u}\in \Lrm^a(0,T;\Lrm^2(\Omega;\R^{d \times m}))$;
\item $u\in \Crm^{0,\gamma}([0,T]\times \overline{\Omega};\R^m)$ for some $\gamma \in (0,1)$;
\item If $p>d$, then $\nabla u\in \Crm^{0,\zeta}([0,T]\times \overline{\Omega};\R^{d \times m})$ for some $\zeta \in (0,1)$.
\end{enumerate}
Moreover, we have the quantitative estimates
\begin{align*}
\norm{\nabla^2 u}_{\Lrm^\infty(\Lrm^p)}&\leq C\Bigl(1+\norm{f}_{\Lrm^\infty(\Lrm^p)}+\norm{f}_{\Lrm^\infty(\Lrm^2)}^{q-1}\Bigr),\\
\norm{\nabla \dot{u}}_{\Lrm^a(\Lrm^2)}&\leq C\Bigl(1+\norm{f}_{\Wrm^{1,a}(\Lrm^2)}\Bigr).
\end{align*}
Here, the constant $C>0$ depends on all constants in the stated assumptions, as well as on $p, a, T, \abs{\Omega}$;
furthermore $\|\cdot\|_{\Lrm^a(\Lrm^p)}$ denotes the norm of the space $\Lrm^a(0,T;\Lrm^p(\Omega))$,
and $\|\cdot\|_{\Wrm^{1,a}(\Lrm^p(\Omega))}$ signifies the norm of the space $\Wrm^{1,a}(0,T;\Lrm^p(\Omega))$.
Denoting by $[u]_{\Crm^{0,\gamma}}$ the $\gamma$-H\"{o}lder seminorm of $u$, the oscillation estimates
stated in (iii) and (iv) above are quantified in the following manner:
\[
[u]_{\Crm^{0,\gamma}([0,T]\times\overline{\Omega})}\leq C\Big(1 + \norm{u}_{\Lrm^\infty(\Wrm^{2,p})}+\norm{\dot{u}}_{\Lrm^a(\Wrm^{1,2})}\Big),
\]
where $\gamma \in (0,1)$, and, if $p>d$,
\[
[\nabla u]_{\Crm^{0,\zeta}([0,T]\times\overline{\Omega})}\leq C\Big(\norm{\nabla u}_{\Lrm^\infty(\Wrm^{1,p})}+\norm{\nabla\dot{u}}_{\Lrm^a(\Lrm^{2})}\Big),
\]
where $\zeta \in (0,1)$. These estimates are explained in Remark~\ref{rem:H} below.
\end{theorem}

\begin{remark}
\label{rem:locex}
If $\Wcal_0$ satisfies~(A\ref{as:funct}) only in a neighborhood of the initial value $u_0$, then a similar existence and regularity result holds, but only with a finite time of existence (in Theorem~\ref{thm:exist} we can show existence and regularity for any $T > 0$ as long as $f$ satisfies (A\ref{as:f}) on $(0,T)$). In fact, the solution can be extended up to the time when $u$ exits the region within which~(A\ref{as:funct}) holds (the \enquote{boundary of convexity}).
\end{remark}
\begin{remark}
\label{rem:scal}
Provided that the tensor $\Abb$ does not depend on the time variable, it can be easily verified that if $u:[0,T]\times \Omega\to \mathbb{R}^m$ is a strong solution to \eqref{eq:PDE}, for given data $u_0,f,\Lcal_t$, then $u_\alpha(t,\cdot):=u(\alpha t,\cdot):[0,\frac{T}{\alpha}]\times \Omega\to \mathbb{R}^m$ is a solution for $u_0,f_\alpha(t,\cdot):=f(\alpha t,\cdot)$ and $\Lcal_{\alpha t}$. We observe
\[
\norm{f_\alpha}_{\Lrm^\infty(0,\frac{T}{\alpha};\Lrm^2(\Omega))}=\norm{f}_{\Lrm^\infty(\Lrm^2)},  \qquad
\norm{f_\alpha}_{\Lrm^a(0,\frac{T}{\alpha};\Lrm^2(\Omega))}= \alpha^{-1/a}\norm{f}_{\Lrm^a(\Lrm^2)},
\]
and
\[
  \norm{\dot{f}_\alpha}_{\Lrm^a(0,\frac{T}{\alpha};\Lrm^2(\Omega))}= \alpha^{1-1/a}\norm{\dot{f}}_{\Lrm^a(\Lrm^2)}.
\]
The first estimate in~(iv) above is invariant under this scaling; the second one can be improved (letting $\alpha \to \infty$
and noting that all constants in the proof of the inequality are independent of $T$) to
\begin{align*}
\norm{\nabla \dot{u}}_{\Lrm^a(\Lrm^2)}&\leq C\norm{\dot{f}}_{\Lrm^a(\Lrm^2)}.
\end{align*}
\end{remark}

We will construct our solutions via a Rothe-type time-discretization scheme and iterative minimization of incremental functionals (see~\eqref{eq:FNk}), as in the theory of Mielke--Theil energetic solutions. In fact, at this time-discrete level the functionals used coincide with the ones in the Mielke--Theil theory. It is only for the corresponding time-continuous limits that differences in the solution concepts are seen to emerge. In connection with this, we refer the reader to the remarks concerning the approximability of rate-independent solutions in~\cite{MielkeRindler09}. Our regularity results follow from \enquote{elliptic} estimates at the discrete level, see Lemma~\ref{lem:apriori-space}. We give a brief formal overview of the relevant estimates in Section~\ref{sc:formal}, which will then be followed by their rigorous proofs in Section~\ref{sc:estimate}.

\begin{theorem} \label{thm:unique}
Under the assumptions~(A\ref{as:first})--(A\ref{as:unique}), the (strong) solution to~\eqref{eq:PDE}, whose existence is guaranteed by Theorem~\ref{thm:exist}, is unique among all strong solutions. Moreover, let $u$ be the solution from Theorem~\ref{thm:exist} and assume that
\[
  v \in \Lrm^1(0,T;(\Wrm^{1,2}_0 \cap \Lrm^q)(\Omega;\R^m))  \qquad\text{with}\qquad
  \dot{v} \in \Lrm^1(0,T;\Wrm^{1,b}_0(\Omega;\R^m)),
\]
where
\[
  \begin{cases}
  b = dp/(dp+p-d)  & \text{if $p \in [2,d)$,} \\
  b > 1            & \text{if $p = d$,} \\
  b = 1            & \text{if $p > d$,}
  \end{cases}
\]
satisfies $v(0)=u_0$ and assume that
\begin{equation} \label{eq:strongsol_weak}
  \int_\Omega \nabla v : \Abb : \nabla \phi+\DD W_0(v)\cdot \phi \dd x
  \leq \int_\Omega R_1(\dot{v}-\phi)-R_1(\dot{v})+f\cdot \phi \dd x,
\end{equation}
for all $\phi \in \Crm^\infty_0(\Omega;\R^m)$. Then, $v = u$.
\end{theorem}

\begin{remark}
The condition~\eqref{eq:strongsol_weak} is in particular satisfied if~\eqref{eq:strongsol_ineq} holds
for solutions $v$ and all $\xi \in \Lrm^1(0,T;(\Wrm^{1,2}_0 \cap \Lrm^q)(\Omega;\R^m))$, see Remark~\ref{rem:exponents}. In the above statement, however, we require even less regularity of $v$ and $\dot{v}$.
\end{remark}

\begin{remark}
Unlike strong solutions, which are unique within their natural regularity class, weak solutions do not have to be unique, as the example in Section~\ref{ssc:ex} demonstrates.
\end{remark}

We will also show under an additional assumption that a continuous piecewise linear finite element approximation of
the model under consideration converges, with a rate, to the unique strong solution. Our main result here is the following theorem.

\begin{theorem}
\label{thm:approx}
Assume~(A\ref{as:first})--(A\ref{as:last}) with $p = 2$ and in the case of $d=3$ additionally that $q \leq 4$.
Then the sequence of numerical solutions $(u_h^\tau)_{h>0, \tau>0}$ generated by the fully discrete
approximation scheme constructed in Section~\ref{sec:num}, on a quasiuniform family of triangulations
of $\Omega$ parametrized by the spatial discretization parameter $h$, and with time step $\tau$, converges to the unique strong solution $u$ of~\eqref{eq:PDE}, which, for some $\tilde\beta \in (0,1)$, satisfies, uniformly in  $h,\tau$,
\[
[u^h_\tau]_{\Crm^{0,\tilde{\beta}}([0,T]\times \Omega)}+ \norm{\nabla u^h_\tau}_{\Lrm^\infty(\Lrm^6(\Omega))}+\norm{\nabla \dot{u}^h_\tau}_{\Lrm^a(\Lrm^2(\Omega))} \leq C.
\]
\end{theorem}

In Section~\ref{sc:rate} we will also prove that the sequence of numerical approximations converges to the unique strong solution $u$ with a rate; see Theorem~\ref{thm:rate}.

\section{Formal a-priori estimates}  \label{sc:formal}

We first illustrate in a non-rigorous way what can be gained from a-priori estimates. These arguments will be made precise in subsequent sections. For the sake of simplicity, we restrict ourselves here to the special case $\Lcal_t = \Delta$. We shall assume therefore that we have a smooth $u \colon [0,T] \times \Omega \to \R^m$ satisfying the initial condition $u(0)=u_0$,
the homogeneous Dirichlet boundary condition $u|_{(0,T] \times \partial \Omega}=0$, and
\begin{equation} \label{eq:toy}
\int_\Omega R_1(\dot{u}(t))-\nabla u(t)\cdot \nabla(\xi-\dot{u}(t))+[-\DD W_0(u(t))+f(t)]\cdot(\xi-\dot{u}(t)) \dd x \leq \int_\Omega R_1(\xi) \dd x
\end{equation}
for all $t \in (0,T]$ and all smooth $\xi \colon \Omega \to \R^m$ such that $\xi|_{(0,T] \times \partial \Omega}=0$.

\subsection{Estimates in space}
\label{ssec:formalspace}
The weak formulation~\eqref{eq:toy} implies that for almost every $(t,x) \in (0,T] \times \Omega$ there is  a $z(t,x)\in \partial R_1(\dot{u}(t,x))$ such that
\[
  z(t) - \Delta u(t)+\DD W_0(u(t))=f(t).
\]
As $R_1$ is globally Lipschitz continuous on $\R^m$ (cf. Assumption (A\ref{as:second})), there exists a positive constant $C$ such that $\abs{z(t,x)}\leq C$, uniformly in $(t,x) \in (0,T] \times \Omega$. Hence,
$- \Delta u(t)\in \Lrm^\infty(\Omega;\R^m)$ if $f$ and $\DD W_0$ are globally bounded.
Moreover, as $\Omega$ is assumed to be bounded with a $\Crm^{1,1}$ boundary, elliptic regularity theory implies that
$\nabla^2 u(t)\in \Lrm^s(\Omega;\R^{d\times d \times m})$ for all $s\in(1,\infty)$, and
\[
\norm{\nabla^2u(t)}_{\Lrm^s}\leq C\norm{\Delta u(t)}_{\Lrm^s} \qquad \mbox{for all $s\in(1,\infty)$},
\]
where the constant $C$ depends only on $\Omega,s,d$; see~\cite[(11.8)]{LadyzhenskayaUraltseva68}, or~\cite{Browder60,Browder61}.
Thus, if for the moment $\DD W_0$ is assumed bounded, we find that
\begin{align*}
 \norm{\nabla^2u(t)}_{\Lrm^s}\leq C(1+\norm{f(t)}_{\Lrm^s}) \qquad \mbox{for all $s\in(1,\infty)$}.
\end{align*}
Morrey's embedding theorem then yields
\[
  \nabla u \in \Lrm^\infty(0,T;\Crm^{0,\alpha}(\cl{\Omega};\R^{d \times m})) \quad \text{for all $\alpha\in[0,1)$.}
\]
The above derivation and the extension to the case when $\DD W_0$ is unbounded but satisfies~\eqref{eq:DWgrowth},~\eqref{eq:Wmon} is made precise in Lemma~\ref{lem:apriori-space}.

\subsection{Testing with $u$}
\label{ssec:u}
Choosing $\xi=\dot{u}(t)-\phi$ for any $\phi\in \Crm^\infty_0(\Omega;\R^m)$ implies that
  \begin{align*}
\int_\Omega R_1(\dot{u}(t))-R_1(\dot{u}(t)-\phi)+\nabla u(t)\cdot \nabla\phi+[\DD W_0(u(t))-f(t)]\cdot\phi \dd x\leq 0.
\end{align*}
Since $R$ is homogeneous of degree 1 and convex, it is subadditive, i.e.\ $R_1(a+b) \leq R_1(a)+R_1(b)$.
Hence,
\begin{equation} \label{eq:toy_phi}
\int_\Omega \nabla u(t)\cdot \nabla\phi+[\DD W_0(u(t))-f(t)]\cdot\phi \dd x\leq \int_\Omega R_1(-\phi)\dd x.
\end{equation}
Using $\phi=u(t)$, and noting
\[
  \DD W_0(u(t)) \cdot u(t) \geq -\mu \abs{u(t)}^2 + \DD W_0(0) \cdot u(t),
\]
which follows from~\eqref{eq:Wmon}, we get that
\begin{align*}
\norm{\nabla u(t)}_{\Lrm^2}^2\leq \mu \norm{u(t)}_{\Lrm^2}^2 + \norm{f(t)}_{\Lrm^2}\norm{u(t)}_{\Lrm^2}+C\norm{u(t)}_{\Lrm^1},
\end{align*}
which implies by Poincar\'{e}'s inequality that
\[
\norm{\nabla u(t)}_{\Lrm^2}\leq \norm{f(t)}_{\Lrm^2}+C(1+\norm{u(t)}_{\Lrm^2})\quad \mbox{for a.e. $t \in (0,T]$}.
\]
When $W_0$ is convex, the term $(\DD W_0(u(t))-\DD W_0(0))\cdot u(t)$ is nonnegative. In this case, or in the case when $\DD W_0$ is globally bounded, we arrive at an upper bound on $\norm{\nabla u(t)}_{\Lrm^2}$ that is independent of $u(t)$.

\subsection{Testing with $\dot{u}$}
An a-priori estimate in a stronger norm can be derived by choosing $\xi=0$. Then,
 \begin{align*}
\int_\Omega R_1(\dot{u}(t))+\partial_t\biggl(\frac{\abs{\nabla u(t)}^2}{2}\biggr)+\partial_t \bigl[ W_0(u(t)) \bigr] \dd x\leq \int_\Omega f(t)\cdot \dot{u}(t)\dd x.
\end{align*}
We integrate this over the time interval $(0,\tau) \subset (0,T)$ and, assuming for simplicity that $u(0) = 0$ and $W_0(0)=0$, we find that
 \begin{align*}
&\int_0^\tau\!\!\int_\Omega R_1(\dot{u})\dd x\dd t+ \int_\Omega\frac{\abs{\nabla u(\tau)}^2}{2}+W_0(\tau)\dd x\\
&\qquad \leq \int_0^\tau\!\!\int_\Omega f\cdot
 \dot{u}\dd x\dd t\\
&\qquad= -\int_0^\tau\!\!\int_\Omega \dot{f}\cdot
 u\dd x\dd t+\int_\Omega f(\tau)\cdot u(\tau) \dd x
\\
&\qquad \leq \norm{\dot{f}}_{\Lrm^1(\Lrm^2)}\norm{u}_{\Lrm^\infty (\Lrm^2)}+\norm{f}_{\Lrm^\infty(\Lrm^2)}\norm{u}_{\Lrm^\infty (\Lrm^2)}.
\end{align*}
This implies, by taking the supremum over all $\tau \in (0,T]$, assuming additionally that $R_1(z)\geq c\abs{z}$, and absorbing via the Poincar\'{e} inequality, that
\begin{align*}
\norm{\dot{u}}_{\Lrm^1(\Lrm^1)}+\norm{\nabla u}^2_{\Lrm^\infty(\Lrm^2)}+\sup_{\tau}\int_\Omega W_0(u(\tau))&\leq C \Bigl( \norm{\dot{f}}_{\Lrm^1(\Lrm^2)}^2+ \norm{f}_{\Lrm^\infty(\Lrm^2)}^2 \Bigr) \dd x\\
&\leq C\norm{f}_{\Wrm^{1,1}(\Lrm^2)}^2.
\end{align*}

\subsection{Testing with $\ddot{u}$}
\label{ssec:formaltime}
The following estimate is the crucial one. Let $t \in (0,T)$ and assume that $\tau \in (0,T-t]$. We will use~\eqref{eq:toy} at $t$ and at $t+\tau$. At time $t$ we take $\xi=\dot{u}(t)+\dot{u}(t+\tau)$ and at time $t+\tau$ we take $\xi=0$. By adding the two inequalities we find
  \begin{align*}
&\int_\Omega R_1(\dot{u}(t+\tau))+R_1(\dot{u}(t))+\nabla (u(t+\tau)-u(t))\cdot \nabla(\dot{u}(t+\tau))\\
&\qquad+[\DD W_0(u(t+\tau))-\DD W_0 (u(t))+f(t)-f(t+\tau)]\cdot \dot{u}(t+\tau) \dd x  \notag\\
&\quad\leq \int_\Omega R_1(\dot{u}(t)+\dot{u}(t+\tau)) \dd x.
\end{align*}
By the subadditivity of $R_1$ we get
  \begin{align*}
&\int_\Omega \frac{\nabla(u(t+\tau)-u(t))}{\tau}\cdot \nabla(\dot{u}(t+\tau))+\biggl[\frac{\DD W_0(u(t+\tau))-\DD W_0 (u(t))}{\tau}\biggr]\cdot\dot{u}(t+\tau) \dd x\\
&\qquad  \leq \int_\Omega \frac{f(t+\tau)-f(t)}{\tau}\cdot\dot{u}(t+\tau) \dd x.
\end{align*}
Letting $\tau\to 0$, and estimating, we find
\begin{align*}
\norm{\nabla \dot{u}(t)}_{\Lrm^2}^2+\int_\Omega \DD^2W_0(u(t))\,[\dot{u}(t),\dot{u}(t)] \dd x \leq \int_\Omega\dot{f}(t)\cdot\dot{u}(t)\dd x
\end{align*}
for any fixed $t \in (0,T)$. We remark that this is exactly the point where the mild convexity
assumption~\eqref{eq:muCP} is essential. Indeed,~\eqref{eq:Wmon} implies (via difference quotients)
in the case where $W_0$ is twice differentiable that
$\DD^2W_0(u)[\dot{u},\dot{u}]\geq -\mu \abs{\dot{u}}^2$.
Hence, by H\"older's inequality and Poincar\'{e}'s inequality,
\begin{align*}
\norm{\nabla \dot{u}(t)}_{\Lrm^2}^2 &\leq \norm{\dot{f}(t)}_{\Lrm^2}\norm{\dot{u}(t)}_{\Lrm^2} +\mu \norm{\dot{u}(t)}_{\Lrm^2}^2\\
&\leq C_P(\Omega)\norm{\dot{f}(t)}_{\Lrm^2}\norm{\nabla\dot{u}(t)}_{\Lrm^2} +\mu C_P(\Omega)^2\norm{\nabla\dot{u}(t)}_{\Lrm^2}^2.
\end{align*}
Thus, invoking~\eqref{eq:muCP} with $\kappa =1$ (recall that we took $\mathcal{L}_t= \Delta$ in this
section, for simplicity), we have
\begin{align*}
\norm{\nabla \dot{u}(t)}_{\Lrm^2}\leq \frac{C_P(\Omega)}{1-\mu C_P(\Omega)^2}\norm{\dot{f}(t)}_{\Lrm^2}.
\end{align*}
Another heuristic approach to the derivation of the above estimate is to use $-(\dot{u}\eta){\dot{~}}$ as test function, where $\eta$ has compact support in $(0,T)$.

\subsection{Optimality in the scalar case}

In the scalar case, corresponding to $m=1$, the following example shows that even if the right-hand side of~\eqref{eq:PDE} is smooth, not more than Lipschitz regularity of the solution can be expected.

\begin{example}
The uniqueness result from Theorem~\ref{thm:unique} implies that, for $t \in [0,2]$, the function $v$ defined by $v(t)=\max\{t-1,0\}$ is the unique strong solution to the rate-independent system
\[
\partial\abs{\dot{v}(t)}\ni t-v(t) \quad\text{for $t \in (0,2]$}\qquad \text{and} \qquad v(0)=0.
\]
Now let
\[
-\Delta \phi +\phi = 1 \quad\text{in $\Omega$} \qquad\text{and}\qquad \phi=0 \quad\text{ on $\partial\Omega$.}
\]
By decomposing $\phi$ as $\phi=\phi_{+} + \phi_{-}$, where $\phi_{+}:=\max\{\phi,0\}$,
$\phi_{-}:=\min\{\phi,0\}$, noting that $\phi_{+}\geq 0$ and $\phi_{-} \leq 0$, and testing the
equation with $\phi_{-}$, we find that $\|\nabla \phi_{-}\|^2_{\Lrm^2} + \|\phi_{-}\|^2_{\Lrm^2} = \int_{\Omega}\phi_{-} \dd x \leq 0$,
whereby $\phi_{-}=0$; hence $\phi=\phi_{+} \geq 0$.
Thus, letting $u(t,x):= v(t) \phi(x)$, we have that
\[
t - u(t) + \Delta u(t) = t-v(t) \in \partial\abs{\dot{v}(t)} \subseteq \partial\abs{\dot{u}(t)}\quad\text{for $t \in (0,2]$}\qquad \text{and} \qquad u(0)=0.
\]
Hence $u$ is a strong solution (in the sense of~\eqref{eq:strongsol_ineq}) to
\[
\partial\abs{\dot{u}}\ni t-u+\Delta u \quad\text{in $(0,2]\times \Omega$} \qquad\text{and}\qquad u(0)=0.
\]
Moreover, clearly $\dot{u}\in \Lrm^\infty([0,2];\Crm^\infty(\overline{\Omega}))\setminus \Crm^0([0,2];\Crm^\infty(\overline{\Omega}))$.
\end{example}

\section{Existence and regularity of solutions}  \label{sc:estimate}

In this section we will prove Theorem~\ref{thm:exist}. We will do so by the Rothe method and discrete analogues of the estimates from Sections~\ref{ssec:formalspace} and \ref{ssec:formaltime}. It will transpire that our a-priori information on solutions is quite strong and thus we obtain compactness in a variety of spaces. The difficulty is to identify the limiting equation, and for this we will need H\"{o}lder continuity of the solution.

\subsection{Time discretization}

We consider a sequence of partitions
\[
  0 = t^N_0 < t^N_1 < \cdots < t^N_N = T,  \qquad\text{where}\qquad
  t^N_k - t^N_{k-1} = \tau:=\frac{T}{N}, \quad N \in \N,
\]
and seek a sequence of approximations
\[
  (u^N_k)_{k=0,\ldots,N} \subset (\Wrm^{1,2}_0 \cap \Lrm^q)(\Omega;\R^m),
\]
which solve a suitable temporally discretized version of~\eqref{eq:PDE}. As a discretization of the source term $f$ we set
\[
  f^N_k := f(t^N_k) \quad\text{for $k = 0,\ldots,N$,}
  \qquad\text{and}\qquad
  f^N := \sum_{k=1}^N \ONE_{(t^N_{k-1},t^N_k]}f^N_k.
\]
Observe that our assumption $f \in \Wrm^{1,a}(0,T;\Lrm^p(\Omega;\R^m))$ for some $a \in (1,\infty)$ implies that
$$f^N\in \Lrm^\infty(0,T;\Lrm^p(\Omega;\R^m)).$$
Further, we define the following approximations of the elliptic operator $\Lcal_t$:
\[
  [\Lcal^N_kv]^\beta = \diverg([\Abb^N_kv]^\beta)
  := \sum_j \partial_j \sum_{i,\alpha} A^{\alpha,\beta}_{i,j}(t_k^N,\frarg)\,\partial_i v^\alpha.
\]

Next, we set $u_{-1}^N:=u_0$ and, successively, for each $k = 0,1,\ldots,N$, minimize the functional
\begin{equation} \label{eq:FNk}
  \Fcal^N_k(v) := \int_\Omega R_1(v-u^N_{k-1}) + \nabla v : \frac{\Abb^N_k}{2} : \nabla v + W_0(v) - f^N_k\cdot v \dd x
\end{equation}
over all $v \in (\Wrm^{1,2}_0 \cap \Lrm^q)(\Omega;\R^m)$. Here, we have used the notation
\[
  \nabla v : \Abb^N_k : \nabla w := \sum_{\alpha,\beta,i,j} A^{\alpha, \beta}_{i,j}(t_k^N,x)\, \partial_i v^\alpha\, \partial_j w^\beta.
\]
First of all, thanks to (A\ref{as:last}), we may set $u_0^N=u_0$. Further, since $R_1$ is convex and lower semicontinuous and $W_0$ is of lower order, we may deduce by the usual Direct Method that a minimizer exists, which we call $u^N_k$. More precisely, we take a minimizing sequence $(v_j)_{j \geq 1} \subset  (\Wrm^{1,2}_0 \cap \Lrm^q)(\Omega;\R^m)$ with $\Fcal^N_k(v_j) \to \min \Fcal^N_k$. Then, by the coercivity of $W_0$ (see~\eqref{eq:Wgrowth}), the strong ellipticity of $\Lcal_t$ (see~\eqref{eq:Lell}) and $R \geq 0$, we get the estimate
\[
  \norm{\nabla v_j}_{\Lrm^2}^2 + \norm{v_j}_{\Lrm^q}^q \leq C(1 + \norm{f^N_k}_{\Lrm^2} \, \norm{v_j}_{\Lrm^2})
\]
for a $j$-independent constant $C > 0$. Thus, using the Poincar\'{e} and Young inequalities,
\[
  \norm{\nabla v_j}_{\Lrm^2} \leq C(1+\norm{f^N_k}_{\Lrm^2})  \qquad\text{and}\qquad
  \norm{v_j}_{\Lrm^q}^q \leq C(1+\norm{f^N_k}_{\Lrm^2}^2).
\]
That is, we have shown coercivity in $(\Wrm^{1,2}_0 \cap \Lrm^q)(\Omega;\R^m)$. Hence, we may assume after selecting a non-relabeled subsequence that $v_j \toweak v$ in $(\Wrm^{1,2} \cap \Lrm^q)(\Omega;\R^m)$. Furthermore, by the compact embedding $\Wrm^{1,2}_0(\Omega;\R^m) \cembed \Lrm^2(\Omega;\R^m)$ and selecting another (not relabeled) subsequence, $v_j \to v$ pointwise almost everywhere. Now, for the convex terms in $\Fcal^N_k$ we get lower semicontinuity immediately, and for $W_0$ ($\geq 0$) we deduce by Fatou's lemma that
\[
  \liminf_{j\to\infty} \int_\Omega W_0(v_j(x)) \dd x
  \geq \int_\Omega W_0(v(x)) \dd x.
\]
Hence, the Direct Method applies and yields the existence of a minimizer, which we call $u^N_k$,
and which satisfies the bounds
\begin{equation} \label{eq:coerc}
  \norm{\nabla u^N_k}_{\Lrm^2} \leq C(1+\norm{f^N_k}_{\Lrm^2})  \qquad\text{and}\qquad
  \norm{u^N_k}_{\Lrm^q}^q \leq C(1+\norm{f^N_k}_{\Lrm^2}^2).
\end{equation}

The minimizer $u^N_k$ satisfies the Euler--Lagrange equation
\[
  0 \in \partial R_1(u^N_k-u^N_{k-1}) - \Lcal^N_k u^N_k + \DD W_0(u^N_k) - f^N_k
\]
in a weak sense. That is, for any test function $\xi \in \Wrm^{1,2}_0(\Omega;\R^m)$ we have that
\begin{align}
  &\int_\Omega R_1(u^N_k-u^N_{k-1}) - \nabla u^N_k : \Abb^N_k : \nabla(\xi-(u^N_k-u^N_{k-1})) \notag\\
  &\qquad + \bigl[ -\DD W_0(u^N_k) + f^N_k \bigr] \cdot (\xi-(u^N_k-u^N_{k-1})) \dd x \leq \int_\Omega R_1(\xi) \dd x.  \label{eq:ELweak}
\end{align}
To see this, we observe that for $\xi \in (\Wrm^{1,2}_0 \cap \Lrm^q)(\Omega;\R^m)$ we have
\begin{equation} \label{eq:Fineq}
  0 \leq \frac{\Fcal^N_k \bigl( u^N_k + \eps(\xi + u^N_{k-1}-u^N_k) \bigr) - \Fcal^N_k(u^N_k)}{\eps},  \qquad \eps > 0.
\end{equation}
First, since $R_1$ is homogeneous of degree 1 and convex, it is subadditive, i.e.\ $R_1(a+b) \leq R_1(a)+R_1(b)$, and so
\begin{align*}
  &R_1 \bigl( u^N_k + \eps(\xi + u^N_{k-1}-u^N_k) - u^N_{k-1} \bigr) - R_1 \bigl(u^N_k - u^N_{k-1} \bigr) \\
  &\qquad = R_1 \bigl( \eps\xi + (1-\eps)(u^N_k - u^N_{k-1}) \bigr) - R_1 \bigl(u^N_k - u^N_{k-1} \bigr) \\
  &\qquad \leq \eps R_1(\xi) - \eps R_1(u^N_k - u^N_{k-1}).
\end{align*}
For the regularizer, using the symmetry of the coefficients in $\Lcal_t$, see~\eqref{eq:Lcoeff}, and setting $\eta := \xi + u^N_{k-1}-u^N_k$, we have that
\begin{align*}
  &\frac{1}{\eps} \int_\Omega [\nabla u^N_k + \eps \nabla \eta] : \frac{\Abb^N_k}{2} : [\nabla u^N_k + \eps \nabla \eta] - \nabla u^N_k : \frac{\Abb^N_k}{2} : \nabla u^N_k \dd x  \\
  &\qquad\to \int_\Omega \nabla u^N_k : \Abb^N_k : \nabla \eta \dd x  \qquad
  \text{as $\eps \todown 0$}
\end{align*}
by the $\Lrm^2$-bounds on all of the quantities involved. Finally, note (using the notation $\dashint_0^\eps:=\frac{1}{\eps}\int_0^\eps$) that
\begin{align*}
  \frac{1}{\eps} \int_\Omega W_0(u^N_k + \eps \eta) - W_0(u^N_k) \dd x
  &= \int_\Omega \dashint_0^\eps \DD W_0(u^N_k + \tau \eta) \cdot \eta \dd \tau \dd x \\
  &\to \int_\Omega \DD W_0(u^N_k) \cdot \eta \dd \tau \dd x
\end{align*}
by the continuity of $\DD W_0$ and the estimate~\eqref{eq:Wgrowth} on the growth of $\DD W_0$. Thus, letting $\eps \todown 0$ in~\eqref{eq:Fineq}, we arrive at~\eqref{eq:ELweak} for $\xi \in (\Wrm^{1,2}_0 \cap \Lrm^q)(\Omega;\R^m)$. A density argument allows us to conclude~\eqref{eq:ELweak} also for $\xi \in \Wrm^{1,2}_0(\Omega;\R^m)$.

Actually, since $u_0^N=u_{-1}^N=u_0$, we find that
\begin{align}
\label{eq:el-u0}
  &\int_\Omega -\nabla u_0 : \Abb_0 : \nabla\xi + \bigl[ -\DD W_0(u_0) + f(0) \bigr] \cdot \xi \dd x \leq \int_\Omega R_1(\xi) \dd x.
\end{align}

\subsection{A-priori estimates in space} \label{ssec:disc-aprior}
For $k=0,\dots,N$, we define the temporal difference quotient
\[
  \delta^N_k:=\frac{u_{k}^N-u_{k-1}^N}{\tau},  \qquad \tau := t^N_k - t^N_{k-1}, \quad t^N_{-1}:=-\tau,
\]
with the convection that $u_0^N=u_{-1}^N=u_0$, as before. Then, upon dividing~\eqref{eq:ELweak} by $\tau$ and replacing $\xi/\tau$ by $\xi$, we get
\begin{align}
&\int_\Omega R_1(\delta_k^N) - \nabla u^N_k : \Abb^N_k : \nabla(\xi-\delta^N_k)+\bigl[ -\DD W_0(u^N_k) + f^N_k \bigr]\cdot(\xi-\delta^N_k) \dd x \notag\\
&\qquad \leq \int_\Omega R_1(\xi) \dd x  \qquad\text{for all $\xi \in \Wrm^{1,2}_0(\Omega;\R^m)$.}  \label{eq:discrete1}
\end{align}
Now, we may further replace $\xi - \delta^N_k$ by $\hat{\xi}$ and use the subadditivity of $R_1$ to get
\[
\int_\Omega - \nabla u^N_k : \Abb^N_k : \nabla \hat{\xi}+\bigl[ -\DD W_0(u^N_k) + f^N_k \bigr] \cdot \hat{\xi} \dd x \leq \int_\Omega R_1(\hat{\xi}) \dd x
\]
for all $\hat{\xi} \in \Wrm^{1,2}_0(\Omega;\R^m)$.

By the 1-homogeneity of $R_1$, we find that $\int_\Omega R_1(\hat{\xi}) \dd x\leq C\norm{\hat{\xi}}_{\Lrm^1}$. Consequently,
\[
\sup_{\{\hat{\xi}\in \Crm^\infty_0(\Omega;\R^m)\,:\,\norm{\hat{\xi}}_{\Lrm^1\leq 1}\}}
\int_\Omega - \nabla u^N_k : \Abb^N_k : \nabla \hat{\xi}+\bigl[ -\DD W_0(u^N_k) + f^N_k \bigr] \cdot \hat{\xi} \dd x \leq C.
\]
Hence, there is a $w^N_k \in \Lrm^\infty(\Omega;\R^m)$, with $\norm{w^N_k}_{\Lrm^\infty} \leq C$ (uniformly in $k,N$), such that
\[
  - \Lcal^N_k u^N_k+\DD W_0(u^N_k) = f^N_k - w^N_k
\]
with the differential operator $\Lcal^N_k$ interpreted in a weak sense.
By Lemma~\ref{lem:apriori-space} below we thus get (recall that $p \in [2,\infty)$ from the assumption on the external force $f$ stated in (A5))
\begin{align}  \label{eq:apri_space}
  \norm{\nabla^2 u^N_k}_{\Lrm^p}\leq C (1+\norm{f^N_k}_{\Lrm^p}+\norm{f^N_k}_{\Lrm^2}^{q-1}).
\end{align}

The next lemma is an elliptic regularity result that is specifically tailored to our situation. We wish to point out that it allows for rather general, physically-motivated, assumptions on the elastic energy functional $\Wcal_0$, namely~\eqref{eq:Wgrowth}--\eqref{eq:Wmon}, but no other structural hypotheses, such as symmetry, are needed.

\begin{lemma} \label{lem:apriori-space}
With our usual hypotheses from Section~\ref{sc:assume}, but excluding the mild convexity assumption~\eqref{eq:muCP}, let $u\in (\Wrm^{1,2}_0 \cap \Lrm^q)(\Omega;\R^m)$ be any weak solution to
\begin{equation} \label{eq:space}
\left\{
\begin{aligned}
-\Lcal_t u +\DD W_0 (u)&=g \quad\text{ in $\Omega$,}\\
u|_{\partial\Omega}&=0.
\end{aligned} \right.
\end{equation}
If $g\in \Lrm^s(\Omega;\R^m)$ for $s\in [2,\infty)$, then
\begin{equation} \label{eq:space_est}
 \norm{\nabla^2u}_{\Lrm^s}\leq C(1+ \norm{g}_{\Lrm^s}+\norm{g}_{\Lrm^{2}}^{q-1}).
\end{equation}
\end{lemma}

Note that the inhomogeneity with exponent $q-1$ is due to the $(q-1)$-growth of $\DD W_0$ via~\eqref{eq:DWgrowth}.

\begin{proof}
The existence of a solution $u\in (\Wrm^{1,2}_0 \cap \Lrm^q)(\Omega;\R^m)$ is guaranteed by the same variational argument as the one above. Analogously to~\eqref{eq:coerc} we find that
\begin{align}
\label{eq:coerc2}
\norm{\nabla u}_{\Lrm^2}\leq C(1+\norm{g}_{\Lrm^2}).
\end{align}

Next, we recall that by regularity theory for elliptic systems, if
\[
-\Lcal_t u =g-\DD W_0 (u)\in \Lrm^s(\Omega;\R^m)
\]
then
\begin{align} \label{eq:reg}
 \norm{\nabla^2u}_{\Lrm^s}\leq C \norm{\Lcal_t u}_{\Lrm^s}
\end{align}
since the coefficients are assumed to be Lipschitz continuous and the boundary of $\Omega$ is of regularity class $\Crm^{1,1}$.
For these results see~\cite{AgmonDouglisNirenberg59,AgmonDouglisNirenberg64} and also~\cite[Theorem 7.3]{GiaquintaMartinazzi12} (the standard scalar case is better known and treated for instance in~\cite{Evans10book}).

Therefore, we are left to establish a bound on $\norm{\DD W_0(u)}_{\Lrm^s}$. In the case $d=2$ we find by Sobolev embedding that $\norm{u}_{\Lrm^s}\leq C\norm{u}_{\Wrm^{1,2}}$ for all $s \in [1,\infty)$. Therefore,
 by~\eqref{eq:DWgrowth},~\eqref{eq:coerc2} and the Poincar\'{e} inequality,
\begin{align}
\label{eq:estw}
  \norm{\DD W_0(u)}_{\Lrm^s} \leq C(1+\norm{u}_{\Lrm^{s(q-1)}}^{q-1})
  \leq C(1+ \norm{g}_{\Lrm^{2}}^{q-1}).
\end{align}
This implies, via~\eqref{eq:reg}, the desired estimate.

In the following we will obtain the same estimate for $d=3$, which we therefore assume from now on until the end of the lemma. We will achieve this goal in several steps.

Our first estimate concerns local $\Wrm^{2,2}$-regularity: for any ball $B_{5R} = B_{5R}(x_0) \subset\Omega$ ($x_0 \in \Omega$, $R > 0$) we will show that
  \begin{align}
  \label{eq:local}
\int_{B_R} \abs{\nabla^2 u}^2 \dd x \leq C \biggl[ \int_{B_{5R}}\frac{\abs{\nabla u}^2}{R^2} + \abs{g}^2\dd x \biggr].
\end{align}
Since the structure of the system~\eqref{eq:space} is invariant under translation and scaling of coordinates, we may assume that $R=1$ and that the ball is centered at the origin, $x_0 = 0$.
Indeed, if $u$ solves~\eqref{eq:space} in the scaled and translated ball $B(z,r)$,
then $\tilde{u}(y)= u(z+Ry)$ solves~\eqref{eq:space} in the unit ball with
$\tilde{g}(y)=R^2g(z+Ry)$ and $\tilde{W}_0:=R^2 W_0$.

Thus, in the following, with no loss of generality, we will consider $u$ to be a solution in the ball $B_5 = B_5(0)$. We take a cut-off function $\eta\in \Crm^1_0(B_2)$, where $B_2 = B_2(0)$, such that $\eta \equiv 1$ on $B_1=B_1(0)$. We define the difference quotient of $v \colon \Omega \to \R^m$ in the direction of the $k$th unit coordinate vector $\ee_k$, $k \in \{1,2,3\}$ and $h \in \R$ by
\begin{align}
\label{eq:dq}
\DD^h_k v(x) := \frac{v(x+h \ee_k)-v(x)}{h}= \dashint_0^h \partial_k v(x+s \ee_k)\dd s.
\end{align}
Now, for $k\in \{1,2,3\}$ and $h \in (0,1/2)$, we take $-\DD^{-h}_k(\eta^2 \DD^h_k(u)) \in \Wrm^{1,2}_0(B_3,\R^m)$ as a test function in~\eqref{eq:space}, with $B_3=B_3(0)$, and employ the summation-by-parts rule for difference quotients, to get
\begin{align*}
(I)+(II) &:= \int_\Omega \sum_{i,j,\alpha,\beta} \DD^h_k (A^{\alpha,\beta}_{i,j}\partial_i u^\alpha)\, \partial_j(\eta^2 \DD^h_ku^\beta)\dd x+ \int_\Omega \DD^h_k(\DD W_0(u))\cdot \DD^h_k u \, \eta^2 \dd x\\
&\phantom{:}= -\int_\Omega g\cdot \DD^{-h}_k (\eta^2 \DD^h_k u)\dd x=:(III).
 \end{align*}
We begin by deriving a lower bound on $(I)$. By the product rule for difference quotients ($\DD^h_k (vw)(x) = \DD^h_k v(x) w(x+h\ee_k) + v(x) \DD^h_k w(x)$) and~\eqref{eq:dq} we find that
\begin{align*}
(I)&\geq   \int_\Omega \sum_{i,j,\alpha,\beta} \DD^h_k(A^{\alpha,\beta}_{i,j})\, \partial_i u^\alpha(\frarg + h\ee_k)\, \partial_j(\eta^2 \DD^h_ku^\beta)\dd x \\
&\qquad +2\int_\Omega \sum_{i,j,\alpha,\beta} A^{\alpha,\beta}_{i,j}\,\DD^h_k(\partial_i u^\alpha)\, \eta \, (\partial_j \eta)\, \DD^h_ku^\beta\dd x\\
&\qquad +\int_\Omega \sum_{i,j,\alpha,\beta} A^{\alpha,\beta}_{i,j}\,\DD^h_k(\partial_i u^\alpha)\, \eta^2 \, \DD^h_k (\partial_j u^\beta) \dd x\\
&\geq - C [\Abb]_{\Crm^{0,1}} \int_\Omega \abs{\nabla u(\frarg + h\ee_k)} \abs{\DD^h_k u} \norm{\nabla \eta}_{\Lrm^\infty} \eta + \abs{\nabla u(\frarg + h\ee_k)}\abs{\DD^h_k\nabla u}\eta^2 \dd x \\
&\qquad -\norm{\nabla \eta}_{\Lrm^\infty}\norm{\Abb}_{\Lrm^\infty}\int_\Omega \abs{\DD^h_k\nabla u}\abs{\DD^h_ku}\eta \dd x + \kappa \int_\Omega \abs{\DD^h_k\nabla u}^2\eta^2 \dd x.
\end{align*}
Young's inequality then implies, recalling that $\mbox{supp } \eta \subset B_2$,
\begin{align}
\label{eq:I}
(I)&\geq
 \frac{\kappa}{2} \int_{B_2} \abs{\DD^h_k\nabla u}^2\eta^2 \dd x - C_{\kappa}(1+\norm{\nabla \eta}_{\Lrm^\infty}^2) [\Abb]_{\Crm^{0,1}}^2\int_{B_3} \abs{\nabla u}^2
+\abs{\DD^h_ku}^2 \dd x.
\end{align}
Next, we estimate $(II)$ by~\eqref{eq:Wmon} to find that
\begin{align}
\label{eq:II}
(II)\geq -\mu\int_{B_2}\abs{\DD^h_k u}^2\eta^2 \dd x.
\end{align}
Finally, choosing $0<\eps<\kappa/2$, we deduce by Young's inequality and~\eqref{eq:dq} that
\begin{align} \label{eq:III}
(III) &\leq C_\eps \int_{B_2} (\abs{g}^2\,\eta^2 +\abs{\DD^{-h}_k \eta}^2 \abs{\DD^{-h}_k u}^2) \dd x \notag\\
&\qquad +\eps\dashint_0^{-h} \int_{B_2} \abs{\DD^h_k \partial_k u(x+s \ee_k)}^2\,\eta^2\dd x \dd s.
 \end{align}
Observe that, as $h\in (0,\frac12)$, we have that
\begin{align*}
\int_{B_2} \abs{\DD^{-h}_k \eta}^2 \abs{\DD^{-h}_k u}^2 \dd x &\leq
\norm{\DD^{-h}_k \eta}^2_{\Lrm^\infty}\,\int_{B_2}\absBB{\dashint_0^{-h} \partial_k u(x+s \ee_k) \dd s}^2 \dd x
\\
&\leq \norm{\partial_k \eta}^2_{\Lrm^\infty}\,
\sup_{s\in (0,h)} \int_{B_2} \abs{\partial_k u(x-s \ee_k)}^2\dd x\\
&\leq
\norm{\nabla \eta}^2_{\Lrm^\infty}\,\int_{B_3} \abs{\nabla u}^2\dd x
\end{align*}
and
\begin{align*}
\dashint_0^{-h} \int_{B_2} \abs{\DD^h_k \partial_k u(x+s \ee_k)}^2\,\eta^2 \dd x \dd s
&\leq \sup_{s\in (0,h)} \int_{B_2} \abs{\DD^h_k \partial_k u(x-s \ee_k)}^2\dd x\\
&\leq \int_{B_3} \abs{\DD^h_k \partial_k u}^2\dd x.
\end{align*}
Combining~\eqref{eq:I},~\eqref{eq:II} and~\eqref{eq:III} with the above, we arrive at
\begin{align*}
\frac{\kappa}{2} \int_{B_1} \abs{\DD^h_k\nabla u}^2 \dd x
&\leq C_{\kappa,\eps}
\int_{B_3}
\abs{\nabla u}^2 \dd x +C_\eps\int_{B_2}\abs{g}^2\dd x + \eps \int_{B_3} \abs{\DD^h_k \nabla u}^2\dd x.
\end{align*}
To conclude the estimate, we need to absorb the $\eps$-term into the left-hand side.
To this end, we note that the previous inequality implies, for $s\in (0,h)$ and $\tau^s_k v(x):=v(x+s \ee_k)-v(x)$, that
 \begin{align*}
 \frac{\kappa}{2} \int_{B_1} \abs{\tau^s_k\nabla u}^2 \dd x
 &\leq C_{\kappa,\eps} h^2
 \int_{B_3}
 \abs{\nabla u}^2 \dd x + C_\eps h^2 \int_{B_2}\abs{g}^2\dd x + \eps \int_{B_3} \abs{\tau^s_k \nabla u}^2\dd x.
\end{align*}
If we transform this back to the translated and scaled ball, we arrive at (now for the original $u$)
 \begin{align*}
 \frac{\kappa}{2} \dashint_{B_R} \abs{\tau^s_k\nabla u}^2 \dd x
 &\leq C_{\kappa,\eps} h^2 \dashint_{B_{3R}} \abs{\nabla u}^2 \dd x +  C_\eps \frac{h^2}{R^2} \dashint_{B_{2R}}\abs{g}^2\dd x + \eps \dashint_{B_{3R}} \abs{\tau^s_k \nabla u}^2\dd x
\end{align*}
for every $R > 0$. Hence we can apply the interpolation result of Giaquinta--Modica type stated in~\cite[Lemma 13]{DieningEttwein08} (with $\gamma_{\nabla u}(R,h) = h^2$ and $\gamma_g(R,h) = h^2R^{-2}$, hence in the result we need to fix $h_0$
as a constant multiple of $R$, e.g.\ $R/10$).

Thus, with a different constant $C_{\kappa,\eps} > 0$, we obtain
\begin{align*}
\int_{B_1} \abs{\DD^h_k\nabla u}^2 \dd x \leq C_{\kappa,\eps} \biggl[ \int_{B_{5}} \abs{\nabla u}^2 + \abs{g}^2\dd x \biggr].
\end{align*}
By letting $h\to 0$ we then deduce the desired local estimate~\eqref{eq:local}.
Moreover, Sobolev embedding then implies that
\begin{align}
\norm{u}_{\Lrm^s(B_R)} &\leq  C\norm{\nabla^2u}_{\Lrm^{2}(B_{R})}+\norm{u}_{\Lrm^1(B_R)} \notag\\
&\leq C \bigg(\frac{1}{R}\norm{\nabla u}_{\Lrm^2(B_{5R})}+\norm{g}_{\Lrm^2(B_{5R})} \bigg)+\norm{u}_{\Lrm^1(B_R)}
\label{eq:local2}
\end{align}
for all $B_{5R}\subset\Omega$ and $s\in (1,\infty)$ (recall that we have assumed that $d=3$).

The aim of the next step is to obtain estimates near the boundary of $\Omega$. More precisely, we will show tangential differentiability up to the boundary, for which we will use a flattening argument. By the same scaling argument as before we assume the center point $0$ to be a boundary point, and that we have a one-to-one diffeomorphism $\Psi \colon B_5 \to \R^3$, such that
\[
  \Psi(B_5\cap [\R^2 \times \{0\}])=\partial\Omega\cap B_5
  \qquad\text{and}\qquad
  \Psi(B^+_5) = \Omega\cap B_5,
\]
where $B^+_R := B_R \cap [\R^2 \times (0,\infty)]$.
By a straightforward transformation, we find that on the half-ball $B^+_5$ we have that $\tilde{u}=u\circ\Psi:B^+_5\to\R^m$ is a weak solution of
\[
\widetilde{\Lcal}_t(\tilde u)+\DD W_0(\tilde u)= g \circ\Psi,
\]
where
 \[
 [\widetilde{\Lcal}_t]^\beta= \sum_j \partial_j\sum_{\alpha,i,k,l} J_{jk} A^{\alpha,\beta}_{kl}J_{li}\,\partial_iu^\alpha,  \qquad J_{ij}:=\partial_i\Psi^j.
 \]

We can assume that $\det \nabla\Psi>\kappa_1$ for some $\kappa_1>0$ depending on the prescribed boundary alone. Thus,~\eqref{eq:Lell} holds with the modulus of ellipticity $\kappa\kappa_1^2$. By~\eqref{eq:local} we get for the half-balls that
   \begin{align}
  \label{eq:local3}
\int_{B_1^+} \abs{\partial_k \nabla \tilde{u}}^2 \dd x  \leq C \biggl[\int_{B_{5}^+}
\abs{\nabla \tilde{u}}^2\dd x+\int_{B_{5}^+}\abs{g \circ \Psi}^2\dd x \biggr]
\end{align}
for $k\in\{1,2\}$. Recalling once again that $d=3$, we will now show that $u\in \Lrm^s(\Omega;\R^m)$ for all $s\in (1,\infty)$ and that \begin{align}
\label{eq:ls}
\norm{u}_{\Lrm^s}\leq C(1+\norm{g}_{\Lrm^2}),
\end{align}
where $C$ depends on $s$, $\Omega$ and the constants of our assumptions.

For $z\in(0,1/2)$ we define
the function
$U^z(x,y):=\int_0^z \partial_z \tilde u(x,y,s)\dd s$.
We first observe that $U^z \in \Wrm^{1,2}(B_{1/2}\cap [\R^2 \times \{0\}])$. Indeed,
\begin{align*}
&\int_{B_{1/2}\cap [\R^2 \times \{0\}]}\absBB{\partial_{x}\int_0^z \partial_z \tilde u(x,y,s)\dd s}^2\dd x\dd y  \\
&\qquad \leq \int_{B_{1/2}\cap [\R^2 \times \{0\}]}\bigg(\int_0^{1/2}\abs{\partial_x\partial_z\tilde u}\dd z\bigg)^2\dd x\dd y.
\end{align*}
 The last estimate holds for $\partial_y$ as well and is in both cases controlled by~\eqref{eq:local3}.
 Sobolev embedding further implies
\[
  \norm{U^z}_{\Lrm^s(B_{1/2}\cap [\R^2 \times \{0\}])}\leq C\norm{U^z}_{\Wrm^{1,2}(B_{1/2}\cap [\R^2 \times \{0\}])} \qquad
  \text{for any $s\in (1,\infty)$.}
\]
From $\tilde u(x,y,z)=\int_0^z \partial_z \tilde u(x,y,s)\dd s$ we therefore deduce that
 \[
 \norm{\tilde{u}}_{\Lrm^s(B_{1/2}^+)}\leq C\norm{\tilde{u}}_{\Wrm^{1,2}(B_{5}^+)} + \norm{g \circ \Psi}_{\Lrm^{2}(B_{5}^+)}.
 \]
To finish, we cover $\Omega$ with finitely many balls. For $x\in \overline{\Omega}$, there exists either $B_R(x)$, such that $B_{5R}(x)\subset \Omega$ or such that $B_{5R}(x)\cap\Omega$ is diffeomorphic to $B^+_{5R}(0)$. Since $\overline{\Omega}$ is compact we can choose a finite subfamily of balls for which either~\eqref{eq:local} or~\eqref{eq:local2} holds.
This enables us to complete the proof for $d=3$ as in case $d=2$ and we get~\eqref{eq:space_est} via~\eqref{eq:reg} and~\eqref{eq:estw}
and a (countable) covering of $\Omega$ with balls $B(y,R)$ such that every $x \in \Omega$ is
in at most $N \in \N$ (with $N$ independent of $x$) balls $B(y,5R)$. In fact, this holds with $N = 11^d$.
This standard covering lemma can be proved by setting $\Omega_0 := \Omega$, choosing
$x_{k+1} \in \Omega_k$ and setting $\Omega_{k+1} := \Omega_k \setminus B(x_k,R)$.
Then, $\abs{x_i - x_j} \geq R$ if $i \neq j$. If $x \in \Omega$ lies in $B(x_k,5R)$ for all
$k \in J \subset \N$, then the reduced balls $B(x_k,R/2)$ are disjoint and all lie in $B(x,5R + R/2)$.
Thus, $J$ can contain at most $11^d$ elements.
\end{proof}

\subsection{Estimates for the time derivatives}
For $k=1,\ldots,N$ we test the $k$th inequality of~\eqref{eq:ELweak} with $\xi=0$ and the $(k-1)$st inequality with $\xi := u^N_k - u^N_{k-2}$. Then we add the two resulting inequalities and divide by $\tau$ to find
\begin{align*}
 &\int_\Omega R_1(\delta_k^N)+R_1(\delta_{k-1}^N)+\nabla u^N_k : \Abb^N_k :  \nabla\delta^N_k-\nabla u^N_{k-1} : \Abb^N_{k-1} :  \nabla\delta^N_k
 \\
 &\qquad+ \bigl[\DD W_0(u^N_k)-\DD W_0(u^N_{k-1})\bigr] \cdot\delta^N_k  -(f^N_k-f^N_{k-1})\cdot\delta^N_k - R_1(\delta^N_k+\delta^N_{k-1}) \dd x \leq 0.
\end{align*}
This can be transformed into
\begin{align*}
 &\int_\Omega R_1(\delta_k^N)+R_1(\delta_{k-1}^N)+(\nabla u^N_k-\nabla u^N_{k-1}) : \Abb^N_k :  \nabla\delta^N_k \\
&\qquad + \bigl[\DD W_0(u^N_k)-\DD W_0(u^N_{k-1})\bigr]\cdot\delta^N_k
-(f^N_k-f^N_{k-1})\cdot\delta^N_k \\
&\qquad- R_1(\delta^N_k+\delta^N_{k-1}) + \nabla u^N_{k-1} : (\Abb^N_k-\Abb^N_{k-1}) :  \nabla\delta^N_k \dd x \quad\leq 0.
\end{align*}
We divide by $\tau > 0$ and use the subadditivity of $R_1$ (to cancel the $R_1$-terms) and the ellipticity~\eqref{eq:Lell} to get
\begin{align*}
 \kappa \int_\Omega \abs{\nabla\delta^N_k}^2 \dd x
  &\leq \int_\Omega \frac{\abs{f^N_k-f^N_{k-1}}}{\tau}\, \abs{\delta^N_k}\dd x+\int_\Omega \frac{\abs{\Abb^N_k-\Abb^N_{k-1}}}{\tau}\, \abs{\nabla u^N_{k-1}}\, \abs{\delta^N_k}\dd x
  \\
  &\qquad- \int_\Omega \frac{1}{\tau^2} \bigl[ \DD W_0(u^N_k) - \DD W_0(u^N_{k-1}) \bigr]  \cdot (u^N_k - u^N_{k-1}) \dd x .
\end{align*}
The first term on the right-hand side can be estimated as follows:
\begin{align*}
  \int_\Omega \frac{\abs{f^N_k-f^N_{k-1}}}{\tau}\,  \abs{\delta^N_k} \dd x
  \leq \biggl( \dashint_{t_{k-1}}^{t_k}\norm{\dot{f}}_{\Lrm^2} \dd t \biggr) \norm{\delta^N_k}_{\Lrm^2} \leq C \biggl( \dashint_{t_{k-1}}^{t_k}\norm{\dot{f}}_{\Lrm^2} \dd t \biggr) \norm{\nabla \delta^N_k}_{\Lrm^2}.
\end{align*}
For the second term we use~\eqref{eq:Lcoeff}, the Poincar\'{e} inequality and~\eqref{eq:coerc} to find
\begin{align}
\label{eq:Acrucual}
\begin{aligned}
\int_\Omega \frac{\abs{\Abb^N_k-\Abb^N_{k-1}}}{\tau}\, \abs{\nabla u^N_{k-1}} \, \abs{\delta^N_k}\dd x
&\leq C [\Abb]_{\Crm^{0,1}([0,T] \times \overline{\Omega})}\norm{\nabla u^N_{k-1}}_{\Lrm^2}\norm{\nabla \delta^N_{k}}_{\Lrm^2} \\
&\leq C (1 + \norm{f^N_{k-1}}_{\Lrm^2})\norm{\nabla \delta^N_{k}}_{\Lrm^2}.
\end{aligned}
\end{align}
For the third term we use~\eqref{eq:Wmon} and the Poincar\'{e} inequality to get
\begin{align*}
  - \frac{1}{\tau^2} \int_\Omega \bigl[ \DD W_0(u^N_k) - \DD W_0(u^N_{k-1}) \bigr]  \cdot (u^N_k - u^N_{k-1}) \dd x
  &\leq \mu \int_\Omega \abs{\delta^N_k}^2 \dd x \\
  &\leq \mu C_P(\Omega)^2 \int_\Omega \abs{\nabla \delta^N_k}^2 \dd x,
\end{align*}
where we recall that $C_P(\Omega) > 0$ denotes the Poincar\'{e} constant of $\Omega$.
Hence, by combining these inequalities, we get
\begin{equation} \label{eq:apri_time}
 \norm{\nabla \delta^N_k}_{\Lrm^2} \leq \frac{C}{\kappa-\mu C_P(\Omega)^2} \biggl[ 1 + \dashint_{t_{k-1}}^{t_k}\norm{\dot{f}}_{\Lrm^2} \dd t +\norm{f^N_k}_{\Lrm^2} \biggr].
\end{equation}
By~\eqref{eq:muCP}, the constant on the right is greater than zero.

\subsection{H\"older continuity of the gradient}
\label{ssec:Hgrad}
In this section only we additionally assume that $p > d$ (see the statement of Theorem~\ref{thm:exist}).
For Borel subsets $E\subset\R^d$ with positive and finite Lebesgue measure we will use the notation
\[
\mean{f}_E:=\dashint_E f \dd x =\frac{1}{\abs{E}}\int_Ef \dd x.
\]
We also define
\[
  u^N(t):=\frac{t-t_{k-1}}{\tau}\,u^N_k+\frac{t_k-t}{\tau}\,u^N_{k-1}  \qquad\text{for $t \in (t_{k-1},t_k]$,\; $k = 0,\ldots,N$,}
\]
and $t_{-1}:=-\tau$.
By taking $k=0$ and recalling that $u^N_0=u^N_{-1}:=u_0$, we have that $u^N(0)=u_0$ for all $N \in \N$.

For $p\in [2,\infty)$, $a\in (1,\infty]$, and $f\in \Wrm^{1,a}(0,T;\Lrm^p(\Omega,\R^m))$ we find by embedding that
$$f\in \Crm^{0,(a-1)/a}([0,T];\Lrm^p(\Omega,\R^m)).$$
Therefore,~\eqref{eq:apri_space} and~\eqref{eq:apri_time} imply that
\begin{align}
\label{eq:apri_uN}
  \norm{\nabla^2 u^N}_{\Lrm^\infty(\Lrm^p)}
  + \norm{\nabla \dot{u}^N}_{\Lrm^a(\Lrm^2)} \leq C,
\end{align}
uniformly in $N$. Thus,
\begin{align}
\label{eq:infty2p}
  u^N \in \Lrm^\infty(0,T;\Wrm^{2,p}(\Omega;\R^m)),  \qquad
  \dot{u}^N \in \Lrm^a(0,T;\Wrm^{1,2}(\Omega;\R^m)),
\end{align}
and their respective norms are uniformly bounded with respect to $N$.

Now, $\Wrm^{2,p}(\Omega;\R^m) \embed \Crm^{1,\alpha}(\cl{\Omega};\R^m)$ for some $\alpha \in (0,1)$ if $p\in (d,\infty)$; hence $u^N\! \in \Lrm^\infty(0,T;\Crm^{1,\alpha}(\cl{\Omega};\R^m))$ and, in fact, the embedding
\[\Lrm^\infty(0,T;\Wrm^{2,p}(\Omega;\R^m)) \cap  \Wrm^{1,a}(0,T;\Wrm^{1,2}(\Omega;\R^m)) \embed \Lrm^r(0,T;\Crm^{1,\alpha}(\cl{\Omega};\R^m)) \]
is compact for any $r \in [1,\infty)$ and any smaller $\alpha \in (0,1)$.
We will show that if $a\in (1,\infty)$ (where the assumption that $a<\infty$ was made with no loss of generality and merely for the sake of simplicity of the exposition) and $p\in (d,\infty)$, then $\nabla u^N$ is uniformly Lipschitz continuous with respect to the metric
\[
  \rho \bigl( (t,x),(s,y) \bigr) := \abs{t-s}^\zeta +\abs{x-y}^\alpha,
\]
for any $\alpha \in (0,1)$ such that $u^N \in \Lrm^\infty(\Crm^{1,\alpha})$, and where
\[
  \zeta = \frac{\alpha(a-1)}{\big(\frac{d}{2}+\alpha\big)a} = \frac{\alpha}{b}, \qquad
  b=\Big(\frac{d}{2}+\alpha\Big)\frac{a}{a-1}.
\]

By Campanato's integral characterization of H\"{o}lder continuity~\cite{Campanato63} (see also Sec.~III.1 in~\cite{Giaquinta83book} or Sec.~2.3 in~\cite{Giusti03book}), we need to show that
\[
  \dashint_{t}^{t+r^b}\!\!\dashint_{B_r(x) \cap \Omega}\abs{\nabla u^N-\mean{\nabla u^N}_{(t,t+r^b)\times B_r(x)}} \dd x\dd t
  \leq Cr^\alpha
\]
for all $(t,t+r^{b})\times B_r(x)\subset [0,T]\times \R^d$, $r > 0$. We remark first that this \enquote{parabolic} version follows from the usual one via the transformation $g(s,x) = bs^{b-1}f(s,x)$. Moreover, from the boundary regularity we infer that $\Omega$ has no internal or external cusps, hence the boundary version of Campanato's Theorem is applicable (see Section~2.3 in~\cite{Giusti03book}). In the following we only prove the statement for internal points, but will comment on the differences for boundary points afterwards.

For ease of exposition we assume that $(t,x)=(0,0)$ and estimate
\begin{align*}
 &\dashint_0^{r^{b}}\!\!\dashint_{B_r}\abs{\nabla u^N-\mean{\nabla u^N}_{(0,r^{b})\times B_r}} \dd x\dd t\\
&\qquad\leq \sup_{0 \leq t \leq r^{b}} \; \dashint_{B_r}\abs{\nabla u^N(t,x)-\mean{\nabla u^N(t)}_{B_r}} \dd x \\
&\qquad\qquad + \dashint_0^{r^{b}}\abs{\mean{\nabla u^N}_{(0,r^{b})\times B_r}-\mean{\nabla u^N(t)}_{B_r}} \dd t\\
&\qquad\leq Cr^\alpha +C \underbrace{\dashint_0^{r^{b}}r^{b}\dashint_{B_r}\abs{\nabla \dot{u}^N} \dd x \dd t}_{=: (I)}.
\end{align*}
Here we have used the a-priori $\Lrm^\infty(\Crm^{1,\alpha})$ regularity result on the first integral and the Poincar\'{e} inequality in the time direction on the second integral.

To bound~$(I)$, we use the $\Lrm^a(\Lrm^2)$ bound on $\nabla \dot{u}^N$ and H\"older's inequality to get with $a' = a/(a-1)$ that
\begin{align*}
 (I)&\leq\int_0^{r^{b}} \Bigg(\dashint_{B_r}\abs{\nabla \dot{u}^N}^2 \dd x\bigg)^{1/2}\dd t \leq r^{b/a'}\bigg(\int_0^{r^{b}}\Bigg(\dashint_{B_r}\abs{\nabla \dot{u}^N}^2 \dd x\bigg)^{a/2}\dd t\bigg)^{1/a}\\
&= Cr^{b/a'-d/2}\norm{\nabla \dot{u}^N}_{\Lrm^a(\Lrm^2)} \leq Cr^\alpha.
\end{align*}
For boundary points $x \in \partial \Omega$, in all of the above estimates $B_r$ needs to be replaced by $B_r(x) \cap \Omega$. In the very last estimate we additionally need to use $\abs{B_r(x) \cap \Omega} \geq \alpha \abs{B_r(x)}$ for some fixed $\alpha \in (0,1)$ (depending on the bounds on the boundary regularity). Thus we may apply Campanato's Theorem to conclude H\"{o}lder regularity.

Finally, Lipschitz continuity in the metric $\rho$ implies H\"{o}lder continuity, with exponent $\zeta$, jointly in space and time. To see this, we note that, since $\zeta< \alpha$,
\begin{align*}
  \abs{t-s}^{\zeta} + \abs{x-y}^\alpha &\leq 2 \max \bigl\{ (\abs{t-s} + \abs{x-y})^{\zeta}, (\abs{t-s} + \abs{x-y})^\alpha \bigr\} \\
  &\leq 2(1+(T + \diam(\Omega))^{\alpha-\zeta}) (\abs{t-s} + \abs{x-y})^{\zeta},
\end{align*}
where we have to consider the cases $\abs{t-s} + \abs{x-y} \leq 1$ and $\abs{t-s} + \abs{x-y} > 1$ separately.

\subsection{H\"older continuity of the solution}
\label{ssec:Hsol}

By a similar argument to the one in the last section, we will show that $u^N$ is uniformly H\"older continuous. We only need to consider the case $p\in [2,d]$, since otherwise the uniform H\"older continuity of $u^N$ follows from~\eqref{eq:infty2p}.

As (see~\eqref{eq:infty2p}) $\nabla^2 u^N \in \Lrm^\infty(0,T;\Lrm^p(\Omega;\R^{m \times d \times d}))$, we have that $\nabla u^N \in \Lrm^\infty(0,T;\Wrm^{1,s}_0(\Omega;\R^{d\times m}))$  for all $s\in [1,\frac{pd}{d-p}]$ when $p<d$, and for all $s\in (1,\infty)$ when $d=p$. Since $\frac{pd}{d-p}\geq 6$, there exists an $\alpha\in (0,1)$ such that $u^N\in \Lrm^\infty(0,T;\Crm^{0,\alpha}(\overline\Omega;\R^m))$. Furthermore, as $\nabla \dot{u}^N\in \Lrm^a(0,T;\Lrm^2(\Omega;\R^{d\times m}))$, we deduce that $\dot{u}^N\in \Lrm^a(0,T;\Lrm^s(\Omega;\R^m))$, for $s\in [1,6]$ if $d=3$ and $s\in (1,\infty)$ for $d=2$. By~\eqref{eq:infty2p},
\[
  u^N \in \Lrm^\infty(0,T;\Crm^{0,\alpha}_0(\overline{\Omega},\R^m)) \cap \Wrm^{1,a}(0,T;\Lrm^{s}(\Omega,\R^m))
\]
and the norms can be correspondingly estimated by a constant independent of $N$.

We can therefore argue exactly as before. Indeed, we will show that $u^N$ is Lipschitz continuous with respect to the metric
\[
  \tilde{\rho} \bigl( (t,x),(s,y) \bigr) := \abs{t-s}^\gamma +\abs{x-y}^\alpha,
\]
for any $\alpha \in (0,1)$ such that $u^N \in \Lrm^\infty(0,T;\Crm^{0,\alpha}(\overline\Omega;\R^m))$ (our $\alpha$ here is different from the one in the previous section), and where
\[
  \gamma = \frac{\alpha(a-1)}{\big(\frac{d}{2}+\alpha\big)a} = \frac{\alpha}{b}, \qquad
  b=\Big(\frac{d}{2}+\alpha\Big)\frac{a}{a-1}.
\]
Again we use the variables $a'=\frac{a}{a-1}$. By Campanato's integral characterization of H\"{o}lder continuity, we need to show that
\[
  \dashint_{t}^{t+r^b}\!\!\dashint_{B_r(x)}\abs{ u^N-\mean{u^N}_{(t,t+r^b)\times B_r(x)}} \dd x\dd t
  \leq Cr^\alpha
\]
for all $(t,t+r^{b})\times B_r(x)\subset [0,T]\times \R^d$, $r > 0$.

We assume that $(t,x)=(0,0)$ and estimate as before
  \begin{align*}
 &\dashint_0^{r^{b}}\!\!\dashint_{B_r}\abs{u^N-\mean{u^N}_{(0,r^{b})\times B_r}} \dd x\dd t\\
&\qquad \leq Cr^\alpha+ C \dashint_0^{r^{b}}r^{b}\dashint_{B_r}\abs{\dot{u}^N} \dd x \dd t \\
&\qquad \leq Cr^\alpha+ C r^\frac{b}{a'}\int_0^{r^{b}}\bigg(\dashint_{B_r}\abs{ \dot{u}^N}^2 \dd x\bigg)^\frac{a}{2} \dd t \\
&\qquad \leq Cr^\alpha.
\end{align*}
Hence, we find that $u^N$ is uniformly $\gamma$-H\"older continuous.
\begin{remark}
\label{rem:H}
The estimates in the last two sub-sections are of general nature. Indeed, what is shown here is that, for $p_0>d$,
\[
[u]_{\Crm^{0,\gamma}([0,T]\times\overline{\Omega})}\leq C\Big(\norm{u}_{\Lrm^\infty(\Wrm^{1,p_0})}+\norm{\dot{u}}_{\Lrm^a(\Lrm^1)}\Big),
\]
for some $\gamma\in (0,1)$.
\end{remark}

\subsection{Proof of Theorem~\ref{thm:exist}}  \label{ssc:proof}

The a-priori estimate~\eqref{eq:apri_uN} implies that there exists a (non-relabeled) subsequence such that
\[
  u^N \toweakstar  u  \quad\text{in $\Lrm^\infty(0,T;\Wrm^{2,p}(\Omega;\R^m))$}
\]
and
\[
  u^N \toweak  u  \quad\text{in $\Wrm^{1,a}(0,T;\Wrm^{1,2}(\Omega;\R^m))$}.
\]
By weak compactness in reflexive Banach spaces we furthermore have
\[
  u^N\toweak u  \quad\text{in $\Lrm^r(0,T;\Wrm^{2,p}(\Omega;\R^m))$}
   \qquad  \text{for any $r\in(1,\infty)$.}
\]

We rewrite~\eqref{eq:discrete1} in a time-continuous form. Observe that on $(t_{k-1},t_k]$ we have $\dot{u}^N=\delta^N_k$ and
\[
 \nabla u^N(t) = \frac{t-t_{k-1}}{\tau}\,\nabla u^N_k+\frac{t_k-t}{\tau}\,\nabla u^N_{k-1}= \nabla u^N_k +\frac{t_k-t}{\tau}\,\nabla (u^N_{k-1}-u^N_{k}).
\]
We also set
\[
  k_N(t) := \setb{k\in\{1,\ldots,N\}}{t\in (t_{k-1},t_k]}  \qquad
  \text{for $t\in [0,T]$.}
\]
Therefore,~\eqref{eq:discrete1} here reads as
\begin{align}
&\int_0^T\!\!\int_\Omega R_1(\dot{u}^N(t)) - \nabla u^N_{k_N(t)} : \Abb^N_{k_N(t)} : \nabla(\xi(t)-\dot{u}^N(t)) \dd x \dd t\notag \\
&\qquad\qquad + \int_0^T\!\!\int_\Omega \bigl[-\DD W_0(u^N_{k_N(t)})+f^N_{k_N(t)} \bigr]\cdot(\xi(t)-\dot{u}^N(t)) \dd x \dd t \notag\\
&\qquad \leq \int_0^T\!\!\int_\Omega R_1(\xi) \dd x\dd t  \label{eq:tolimit}
\end{align}
for all $\xi \in \Lrm^1(0,T;\Wrm^{1,2}_0(\Omega;\R^m))$ (first use only $\xi$ that are piecewise constant as a function of $t$ with respect to the subdivision $\{t^N_0,t^N_1,\ldots,t^N_N\}$, and then argue by density).

Using the H\"older continuity of $u$, we find by the Arzel\`{a}--Ascoli theorem a subsequence such that
\[
  u^N \to u
  \quad\text{in $\Crm^{0,\gamma}([0,T]\times\overline{\Omega};\R^m)$}\qquad
  \text{for $0<\gamma< \frac{\alpha(a-1)}{\big(\frac{d}{2}+\alpha\big)a}$}\qquad
  \text{if $p\in [2,d]$}
\]
and
\[
  u^N \to u
  \quad\text{in $\Crm^{1,\zeta}([0,T]\times\overline{\Omega};\R^m)$}\qquad
  \text{for $0<\zeta< \frac{\alpha(a-1)}{\big(\frac{d}{2}+\alpha\big)a}$}\qquad
  \text{if $p\in (d,\infty]$.}
\]
Here, $\alpha$ is defined as above via Morrey's embedding theorem with respect to the spatial variables; see Remark~\ref{rem:H}. By equicontinuity we also know that in both cases
\[
  u^N_{k_N(t)}\to  u(t,\frarg)  \quad\text{in $\Crm^{0,\gamma}(\overline{\Omega};\R^m)$} \quad
  \text{for every $t \in (0,T]$ and any $0<\gamma<\frac{\alpha(a-1)}{\big(\frac{d}{2}+\alpha\big)a}$.}
\]

By the convexity and lower semicontinuity of $R_1$ as well as the assumptions on $\DD W_0$ (continuity) and $f \in \Crm^0(0,T;\Lrm^2(\Omega;\R^m))$, we get
\begin{align*}
&\int_0^T\!\!\int_\Omega R_1(\dot{u})+\bigl[-\DD W_0(u)+f\bigr]\cdot(\xi-\dot{u}) \dd x \dd t\\
&\qquad \leq \liminf_{N\to\infty}\int_0^T\!\!\int_\Omega R_1(\dot{u}^N)+ \bigl[-\DD W_0(u^N_{k_N})+f^N_{k_N} \bigr]\cdot(\xi-\dot{u}^N) \dd x \dd t
\end{align*}
for all $\xi\in \Lrm^1(0,T;\Wrm^{1,2}_0(\Omega;\R^m))$.

The term of the regularizer needs special attention. Rellich's compactness theorem implies that $\Lrm^\infty(0,T;\Wrm^{2,p}(\Omega;\R^m)) \cap \Wrm^{1,a}(0,T;\Wrm^{1,2}_0(\Omega;\R^m))$ is compactly embedded in $\Crm^{0,\beta}([0,T];\Wrm^{1,2}(\Omega;\R^m))$ for some $\beta>0$; see~\cite{Simon87}. Therefore by passing to yet another (not relabeled) subsequence, we find that $\nabla u^N\to\nabla u$ in the strong topology of $\Crm^{0,\beta}([0,T];\Wrm^{1,2}(\Omega;\R^m))$, in particular
\[
  \norm{\nabla u^N(t) - \nabla u(t)}_{\Lrm^2} \to 0 \quad\text{uniformly in $t \in [0,T]$.}
\]
Consequently,
\[
\int_0^T\!\!\int_\Omega\nabla u^N_{k_N} : \Abb^N_{k_N} : \nabla (\xi-\dot{u}^N)\dd x\dd t\to \int_0^T\!\!\int_\Omega\nabla u: \Abb : \nabla(\xi-\dot{u})\dd x\dd t.
\]
Hence, letting $N\to\infty$ in~\eqref{eq:tolimit}, we get
\begin{align*}
&\int_0^T\!\!\int_\Omega R_1(\dot{u}) - \nabla u : \Abb : \nabla(\xi-\dot{u})+\bigl[-\DD W_0(u)+f\bigr]\cdot(\xi-\dot{u}) \dd x\dd t \\
&\qquad \leq \int_0^T\!\!\int_\Omega R_1(\xi) \dd x\dd t,
\end{align*}
for all $\xi \in \Lrm^1(0,T;\Wrm^{1,2}_0(\Omega;\R^m))$. Therefore, the limit inequality~\eqref{eq:strongsol_ineq} is established and our $u$ is indeed a strong solution to~\eqref{eq:PDE}.

We are left to show that the initial value is attained.
Since we already know that $u$ is uniformly continuous on $[0,T]\times \overline{\Omega}$, it remains to show that $u$ actually attains the initial value $u_0$. To see this, we
note that all members of the sequence $(u^N(t))_{N \geq 1}$ attain the initial datum $u_0$ at $t=0$. Since the sequence $(u^N)_{N \geq 1}$ converges uniformly to $u$ in $\Crm^{0,\gamma}([0,T]\times\overline{\Omega};\R^m)$, the initial value $u_0$ is attained by $u$ as well. That completes the proof of Theorem~\ref{thm:exist}.

\section{Uniqueness}

Let $v$ be as in the statement of Theorem~\ref{thm:unique}. One can check that our regularity assumptions on $v$ stated there ensure that all of the quantities appearing below are well-defined.
To prove Theorem~\ref{thm:unique}, we proceed as follows. First, we prove additional regularity for $v$ (for $u$ from Theorem~\ref{thm:exist} this is already known). Then, equipped with this improved regularity, we can use a Gronwall argument to show that $u = v$.

\subsection{Additional regularity}

We take $\epsilon\phi$ with $\phi\in \Crm^\infty_0(\Omega;\R^m)$ as a test function in~\eqref{eq:strongsol_weak} and divide by $\epsilon>0$ to get
 \begin{align*}
 \int_\Omega \nabla v : \Abb : \nabla \phi+\DD W_0(v)\cdot \phi \dd x&\leq \int_\Omega \frac{R_1(\dot{v}-\epsilon \phi)-R_1(\dot{v})}{\epsilon}+f\cdot \phi \dd x
 \\
 &=\int_\Omega R\Big(\frac{\dot{v}}{\epsilon}- \phi\Big)-R\Big(\frac{\dot{v}}{\epsilon}\Big)+f\cdot \phi \dd x
\\
 &\leq C(1+\norm{f}_{\Lrm^\infty(\Lrm^p)}) \norm{\phi}_{\Lrm^{p/(p-1)}},
 \end{align*}
 where for the last estimate we also used $R_1(a+b)-R_1(a)\leq R_1(b)$, which follows from the $1$-homogeneity and convexity of $R_1$. By the dual formula for norms, we therefore find
 \[
 \norm{-\Lcal v +\DD W_0(v)}_{\Lrm^\infty(\Lrm^p)}\leq C(1+\norm{f}_{\Lrm^\infty(\Lrm^p)}).
 \]
Lemma~\ref{lem:apriori-space} then implies that $\nabla^2 v \in \Lrm^\infty(0,T;\Lrm^p(\Omega;\R^m))$ and thus, by Sobolev embedding,
\[
  v \in \Lrm^\infty(0,T;\Lrm^\infty(\Omega;\R^{d\times m})), \qquad
  \nabla v \in \Lrm^\infty(0,T;\Lrm^{b'}(\Omega;\R^{d\times m})),
\]
where $b'$ is the conjugate exponent to $b$, as in the statement of Theorem~\ref{thm:unique} (one can check this case-by-case).

This shows in particular that $\dpr{\nabla v, \nabla \dot{u} - \nabla \dot{v}}$ and $\dpr{\DD W_0(v), \dot{u} - \dot{v}}$ are in $\Lrm^1(0,T)$; the same holds if the roles of $u$ and $v$ are interchanged. Thus, also $v$ satisfies the inequality for strong solutions~\eqref{eq:strongsol_ineq} (even though the regularity of $v$ is different from the one required for $u$).

Finally, since we have assumed that $\nabla \dot{v} \in \Lrm^1(0,T;\Lrm^b(\Omega;\R^{d\times m}))$ and $b \geq 1$, it then also follows that $\dot{v} \in \Lrm^1(0,T;\Lrm^{d/(d-1)}(\Omega;\R^m))$, which for $d \in \{2,3\}$ yields that
\begin{equation} \label{eq:dotv_L32}
  \dot{v} \in \Lrm^1(0,T;\Lrm^{3/2}(\Omega;\R^m)).
\end{equation}

\subsection{A Gronwall argument}

We test~\eqref{eq:strongsol_ineq} for $u$ with $\xi := \dot{v}$ and~\eqref{eq:strongsol_ineq} for $v$ with $\xi := \dot{u}$, integrate from $0$ to $\tau \in (0,T]$ in time and add the resulting inequalities. This gives
\begin{align*}
&\int_0^\tau\int_\Omega R_1(\dot{u})+R_1(\dot{v})- \nabla u : \Abb : \nabla(\dot{v}-\dot{u}) -\nabla v : \Abb : \nabla(\dot{u}-\dot{v})  \dd x\dd t  \\
& \qquad\qquad + \int_0^\tau\int_\Omega \bigl[\DD W_0(v)-\DD W_0(u)\bigr]\cdot(\dot{v}-\dot{u}) \dd x\dd t \\
&\qquad \leq \int_0^\tau\int_\Omega R_1(\dot{u})+R_1(\dot{v}) \dd x\dd t.
\end{align*}
Consequently,
\begin{align}
&\int_0^\tau\frac{\di}{\di t} \int_\Omega \frac{\nabla (v-u) : \Abb : \nabla(v-u)}{2}\dd x\dd t \\
&\qquad\qquad + \int_0^\tau \int_\Omega \bigl[\DD W_0(v)-\DD W_0(u)\bigr]\cdot(\dot{v}-\dot{u})\dd x \dd t  \notag \\
&\qquad\leq \int_0^\tau\int_\Omega \nabla (v-u) : \frac{\partial_t\Abb}{2}: \nabla(v-u)\dd x \dd t  \notag \\
&\qquad \leq C \int_0^\tau\norm{\nabla (u-v)}_{\Lrm^2}^2 \dd t.  \label{eq:gw0}
\end{align}
We wish to use Gronwall's lemma. Hence, we estimate further by the conditions~\eqref{eq:Wmon} and the symmetry of $\DD^2 W_0$, to find
\begin{align}
 &\bigl[\DD W_0(v)-\DD W_0(u)\bigr]\cdot(\dot{v}-\dot{u})  \notag \\
 &\qquad =\int_0^1\DD^2 W_0(\theta u+(1-\theta)v) [u-v,\dot{u}-\dot{v}] \dd \theta  \notag\\
 &\qquad =\frac{\di}{\di t}\bigg(\frac{1}{2}\int_0^1\DD^2 W_0(\theta u+(1-\theta)v)[u-v,u-v] \dd \theta\bigg)   \notag\\
 &\qquad\qquad -\int_0^1 \partial_t\big( \DD^2W_0(\theta u+(1-\theta)v) \bigr) [u-v,u-v] \dd \theta . \label{eq:gw1}
\end{align}
By~\eqref{eq:DDDW} a possible singularity at $0$ in the innermost integral is integrable. Thus, plugging~\eqref{eq:gw1} into~\eqref{eq:gw0},
and also using the regularity assumptions, H\"older's inequality and the Sobolev--Poincar\'{e} inequality, we estimate
\begin{align*}
(I) &:=\frac{1}{2}\int_0^\tau \frac{\di}{\di t} \int_\Omega \biggl( \nabla (v-u) : \Abb : \nabla(v-u) \\
&\qquad\qquad\qquad\qquad + \int_0^1\DD^2 W_0(\theta u+(1-\theta)v)\,[u-v,u-v] \dd \theta \biggr) \dd x \dd t \\
&\phantom{:} \leq \int_0^\tau\int_\Omega\absBB{\int_0^1 \partial_t \bigl( \DD^2 W_0(\theta u+(1-\theta)v) \bigr) \dd \theta} \, \abs{u-v}^2 \dd x \dd t \\
&\qquad+ C \int_0^\tau\norm{\nabla (u-v)}_{\Lrm^2}^2\dd t\\
&\phantom{:} \leq \int_0^\tau\int_\Omega\int_0^1 \abs{\DD^3 W_0(\theta u+(1-\theta)v)}\dd \theta \, (\abs{\dot{v}}+\abs{\dot{u}})\,\abs{u-v}^2\dd x\dd t\\
&\qquad + C \int_0^\tau\norm{\nabla (u-v)}_{\Lrm^2}^2\dd t \\
&\phantom{:} \leq C \int_0^\tau(\norm{\dot{u}}_{\Lrm^\frac{3}{2}}+\norm{\dot{v}}_{\Lrm^\frac{3}{2}})\,\norm{u-v}_{\Lrm^6}^2\dd t+C \int_0^\tau\norm{\nabla (u-v)}_{\Lrm^2}^2\dd t \\
&\phantom{:} \leq C \int_0^\tau (1+\norm{\dot{u}}_{\Lrm^\frac{3}{2}}+\norm{\dot{v}}_{\Lrm^\frac{3}{2}}) \,\norm{\nabla(u-v)}_{\Lrm^2}^2\dd t.
\end{align*}
On the other hand, using the fact that $(u-v)(0)\equiv0$, we find by~\eqref{eq:Lell},~\eqref{eq:muCP} that
\begin{align}
(I)&=\frac{1}{2} \biggl[ \int_{\Omega} \nabla (v-u) : \Abb : \nabla(v-u) \dd x  \notag\\
&\qquad\qquad+ \int_{\Omega} \int_0^1\DD^2 W_0(\theta u+(1-\theta)v)\,[u-v,u-v] \dd \theta\dd x \biggr]_{t=\tau}  \notag\\
&\geq \biggl[ \frac{1}{2}\int_{\Omega} \kappa\abs{\nabla (u-v)}^2-\mu\abs{u-v}^2\dd x \biggr]_{t=\tau}  \notag\\
&=(\kappa-\mu C_p(\Omega)^2) \, \norm{\nabla (u(\tau)-v(\tau))}_{\Lrm^2}^2. \label{eq:gw3}
\end{align}
Therefore, defining
\[
  g := \norm{\nabla (u-v)}_{\Lrm^2}^2 \in \Lrm^\infty(0,T)  \qquad\text{and}\qquad
  h := 1+\norm{\dot{u}}_{\Lrm^\frac{3}{2}}+\norm{\dot{v}}_{\Lrm^\frac{3}{2}} \in \Lrm^1(0,T),
\]
where the integrability follows from~\eqref{eq:dotv_L32}, we find that
\[
g(\tau)\leq c \int_0^\tau h(t) g(t) \dd t \qquad\text{and}\qquad
 g(0)=0,
\]
which implies that $g\equiv 0$ by Gronwall's lemma. Thus, $u \equiv v$ since $(u-v)(t)$ has zero trace on $\partial\Omega$
for all $t \in (0,T]$. That completes the proof of Theorem~\ref{thm:unique}. \qed

\section{Finite element approximation}
\label{sec:num}

The aim of this section is to introduce an approximation scheme for~\eqref{eq:PDE}, based on a spatial finite element discretization of the problem under consideration.

\subsection{Construction}
We shall consider the family of finite-dimensional spaces $(\Xrm^h)_{h>0}$, with $\Xrm^h \subset
\Crm^{0,1}_0(\overline{\Omega};\R^m)$,
parametrized by the discretization parameter $h \in (0,1]$; $\Xrm^h$ will be supposed to be a finite element space generated by $m$-component continuous piecewise
polynomial vector functions on a quasiuniform partition $\mathcal{T}^h$ of $\overline\Omega$ into closed $d$-dimensional simplices $K$ (with, possibly,
curved faces for simplices $K$ that have nonempty intersection with $\partial\Omega$), with $h=\max_{K \in \mathcal{T}^h} h_K$, where $h_K:=\mbox{diam } K$,
and such that $\bigcup_{0<h\leq 1} \Xrm^{h}$ is dense in $\Wrm^{1,2}_0(\Omega;\R^m)$.

We adopt the following standard assumption on $\Xrm^h$: we assume the existence of a projector $\Pbb^h: \Wrm^{1+s,p}_0(\Omega;\R^m)\to \Xrm^h$,
that satisfies, uniformly in $h$, the following approximation result:
\begin{alignat}{2}
\label{eq:proj1}
\|g - \Pbb^h g\|_{\Wrm^{\sigma,p}} \leq C h^{1+s-\sigma}\|g\|_{\Wrm^{1+s,p}},
\end{alignat}
where $0 \leq \sigma \leq 1+s,\quad 0 \leq s \leq 1, \quad 1 < p < \infty$.

As before, we consider a sequence of partitions of the time interval $[0,T]$:
\[
  0 = t_0 < t_1 < \cdots < t_N = T,  \qquad\text{where}\qquad
  t_k - t_{k-1} = \frac{T}{N}=:\tau, \quad N \in \N.
\]
Our time stepping scheme will be based on the implicit Euler method; specifically, we aim to find
\[
  (u^h_k)_{k=0,\ldots,N} \subset \Xrm^h,
\]
which solves a suitable discrete version of~\eqref{eq:PDE}. Let us introduce the following notation for the values of the
external force $f$ at the temporal mesh points:
\[
  f_k := f(t_k) \quad\text{for $k = 0,\ldots,N$,}
\]

Observe that our assumption $f \in \Wrm^{1,a}(0,T;\Lrm^p(\Omega;\R^m))$ implies that $f_k$ is well defined.
To ensure compatibility of the initial datum we define the following elliptic projector:
\[
\Pi^h_0: \Wrm^{1,2}_0(\Omega;\R^m)\to \Xrm^h\quad\text{ defined as }\quad\skp{\Abb_{0}\nabla (g-\Pi^h_0(g))}{\nabla \phi}=0\quad \text{ for all } \phi\in \Xrm^h.
\]
We then let
\begin{align}
\label{eq:uh0}
u^h_{-1} = u^h_0 :=\Pi^h_0(u_0).
\end{align}

Now, successively at each $k = 0,1,\ldots,N$, we minimize the functional
\[
  \Fcal^h_k(v^h) := \int_\Omega R_1(v^h-u^h_{k-1}) + \nabla v^h : \frac{\Abb_{k}}{2} : \nabla v^h + W_0(v^h) - f_k\cdot v^h \dd x
\]
over all $v^h \in \Xrm^h$. Here, we have used the notation
\[
  \nabla v^h : \Abb_{k} : \nabla w^h := \sum_{\alpha,\beta,i,j} A^{\alpha, \beta}_{i,j}(t_k,x)\, \partial_i (v^h)^\alpha \,\partial_j (w^h)^\beta.
\]
Since $R_1$ is convex and lower semicontinuous and $W_0$ is of lower order, we may deduce by the usual Direct Method that a minimizer exists, which we call $u^h_k$. More precisely, the positivity of $R_1$, the ellipticity of $\Lcal$ and the coercivity of $W_0$ imply that a minimizing sequence $v_j^h\subset \Xrm^h$ satisfies the following uniform bound:
\[
  \norm{\nabla v_j^h}_{\Lrm^2}^2 + \norm{v_j^h}_{\Lrm^q}^q \leq C(1 + \norm{f_k}_{\Lrm^2} \, \norm{v_j^h}_{\Lrm^2})
\]
with a $j$-independent constant $C > 0$. Thus, using the Poincar\'{e} and Young inequalities,
\begin{equation*}
  \norm{\nabla v_j^h}_{\Lrm^2} \leq C(1+\norm{f_k}_{\Lrm^2})  \qquad\text{and}\qquad
  \norm{v_j^h}_{\Lrm^q}^q \leq C(1+\norm{f_k}_{\Lrm^2}^2).
\end{equation*}
 Hence, as $\Xrm^h$ has finite dimension, we deduce after selecting a subsequence of the sequence $(v_j^h)_{j \geq 1}$ (not indexed) that $v_j^h \to v^h$ in $\Xrm^h$, as well as pointwise, and $v^h=: u^h_{k}$ is a minimizer that satisfies
 \begin{equation} \label{eq:h-coerc-d}
 \norm{u^h_{k}}_{\Lrm^q}^q + \norm{\nabla u^h_{k}}_{\Lrm^2}^2 \leq C(1+\norm{f(t_k)}_{\Lrm^2}^2).
\end{equation}

Next we show, that $u^h_k$ satisfies the following Euler--Lagrange equation
\[
  0 \in \partial R_1(u^h_k-u^h_{k-1}) - \Lcal^h_k u^h_k + \DD W_0(u^h_k) - f_k,
\]
in a suitable sense. That is, for any test function $\xi^h \in \Xrm^h$ we have that
\begin{align}
  &\int_\Omega R_1(u^h_k-u^h_{k-1}) - \nabla u^h_k : \Abb_{k} : \nabla(\xi^h-(u^h_k-u^h_{k-1})) \notag\\
  &\qquad + \bigl[ -\DD W_0(u^h_k) + f_k \bigr] \cdot (\xi^h-(u^h_k-u^h_{k-1})) \dd x \leq \int_\Omega R_1(\xi^h) \dd x.  \label{eq:ELweak-d}
\end{align}
To see this, we observe that first for $\xi^h \in \Xrm^h$ we have
\begin{equation} \label{eq:Fineq-d}
  0 \leq \frac{\Fcal^h_k \bigl( u^h_k + \eps(\xi^h + u^h_{k-1}-u^h_k) \bigr) - \Fcal^h_k(u^h_k)}{\eps},  \qquad \eps > 0.
\end{equation}
Next, since $R_1$ is homogeneous of degree 1 and convex, it is subadditive, i.e.\ $R_1(a+b) \leq R_1(a)+R_1(b)$, and so
\begin{align*}
  &\frac{1}{\epsilon}\Big(R \bigl( u^h_k + \eps(\xi^h + u^h_{k-1}-u^h_k) - u^h_{k-1} \bigr) - R \bigl(u^h_k - u^h_{k-1} \bigr) \Big)\\
  &\qquad = \frac1{\epsilon} \Big(R \bigl( \eps\xi^h + (1-\eps)(u^h_k - u^h_{k-1}) \bigr) - R \bigl(u^h_k - u^h_{k-1} \bigr) \Big)\\
  &\qquad \leq  R_1(\xi^h) -  R_1(u^h_k - u^h_{k-1}).
\end{align*}
For the regularizer we may compute using the symmetry of the coefficients in $\Lcal_t$, see~\eqref{eq:Lcoeff}, and setting $\eta^h := \xi^h + u^h_{k-1}-u^h_k$,
\begin{align*}
  &\frac{1}{\eps} \int_\Omega [\nabla u^h_k + \eps \nabla \eta^h] : \frac{\Abb_{k}}{2} : [\nabla u^h_k + \eps \nabla \eta^h] - \nabla u^h_k : \frac{\Abb_{k}}{2} : \nabla u^h_k \dd x
  \\
  &\qquad
  \to \int_\Omega \nabla u^h_k : \Abb_{k} : \nabla \eta^h \dd x  \qquad
  \text{as $\eps \todown 0$}.
\end{align*}
 Finally, we note that
\begin{align*}
  &\frac{1}{\eps} \int_\Omega W_0(u^h_k + \eps \eta^h)-f_k\cdot(u^h_k+\eps\eta^h) - W_0(u^h_k)+f_k\cdot u^h_k\dd x
  \\
  &\qquad = \int_\Omega \dashint_0^\eps \DD W_0(u^h_k + s \eta^h) \cdot \eta \dd s\dd x +\skp{f_k}{\eta^h} \\
 &\qquad\to \int_\Omega \DD W_0(u^h_k) \cdot \eta^h \dd \tau \dd x+\skp{f_k}{\eta^h}
\end{align*}
by the continuity of $\DD W_0$. Thus, letting $\eps \todown 0$ in~\eqref{eq:Fineq-d}, we arrive at~\eqref{eq:ELweak-d} for all $\xi^h \in \Xrm^h$.

\subsection{Discrete a-priori estimates in space} \label{ssec:disc-aprior-d}
For $k=0,\ldots,N$ set
\[
  \delta^{\tau,h}_k:=\frac{u_{k}^h-u_{k-1}^h}{\tau},  \qquad \tau := \frac{T}{N}.
\]
Then, dividing~\eqref{eq:ELweak-d} by $\tau$ and replacing $\xi^h/\tau$ by $\xi^h$, we get
\begin{align}
&\int_\Omega R_1(\delta^{\tau,h}_k) + \nabla u^h_k : \Abb_{k} : \nabla(\delta^{\tau,h}_k-\xi^h)+\bigl[-\DD W_0(u^h_k) + f_k \bigr]\cdot(\xi^h-\delta^{\tau,h}_k) \dd x \notag\\
&\qquad \leq \int_\Omega R_1(\xi^h) \dd x.  \label{eq:discrete1-d}
\end{align}
Now, we may further replace $ \xi^h-\delta^{\tau,h}_k$ by $\hat{\xi}^h$ and use the subadditivity of $R_1$ at $\xi = \hat{\xi} + \delta^N_k$ to get
\[
- \int_\Omega  \nabla u^h_k : \Abb_{k} : \nabla \hat{\xi}^h+\bigl[ -\DD W_0(u^h_k) + f_k \bigr] \cdot \hat{\xi}^h \dd x \leq \int_\Omega R_1(\hat{\xi}^h) \dd x.
\]

Given any $z^h \in \Xrm^h$, we define $\mathcal{L}_k^h z^h \in \Xrm^h$ as the (unique) solution of the problem
\begin{equation}\label{eq:Lkh}
\skp{\mathcal{L}_k^h z^h}{\hat{\xi}^h}=-\int_\Omega \nabla z^h :\Abb_{k}:\nabla \hat{\xi}^h\, \dd x \qquad \text{for all $\hat{\xi}^h \in \Xrm^h$.}
\end{equation}
We then define $\psi^h_k:=\mathcal{L}_k^hu_k^h$. Using this test function we have that
\begin{align}
\label{eq:testlaplace}
\skp{\psi^h_k}{\psi^h_k} + \int_\Omega \bigl[- \DD W_0(u^h_k) + f_k \bigr] \cdot \psi^h_k \dd x \leq \int_\Omega R_1(\psi_h^k) \dd x.
\end{align}
By H\"older's inequality we deduce from~\eqref{eq:testlaplace} that
\begin{align*}
\norm{\psi^h_k}_{\Lrm^2}\leq \norm{\DD W_0(u^h_k)}_{\Lrm^2}+\norm{f_k}_{\Lrm^2}+C.
\end{align*}
By~\eqref{eq:DWgrowth} and Sobolev embedding, assuming additionally in the case of $d=3$ that $q \leq 4$, we find
that
\begin{align*}
\int_{\Omega}\abs{\DD W_0(u^h_k)}^2\, \dd x&\leq C \int_{\Omega} 1+\abs{u^h_k}^{2(q-1)}\, \dd x
\leq C\bigg(\int_{\Omega}\abs{\nabla u^h_k}^{2}\, \dd x\bigg)^{q-1}+C.
\end{align*}
Therefore, by~\eqref{eq:h-coerc-d}
\begin{align*}
\norm{\psi^h_k}_{\Lrm^2}\leq C\,\norm{f_k}_{\Lrm^2}^{\max\{1,(q-1)/2\}}+C\leq C\,\norm{f}_{\Lrm^\infty(\Lrm^2)}^{\max\{1,(q-1)/2\}}+C.
\end{align*}
Hence,
\begin{align}
\label{eq:laplacelp}
\norm{-\mathcal{L}_k^hu_k^h}_{\Lrm^2} \leq C\,\norm{f}_{\Lrm^\infty(\Lrm^2)}^{\max\{1,(q-1)/2\}}+C.
\end{align}
Now we note that (cf.\ the Appendix for a proof):
\begin{alignat}{2}
\begin{aligned}
\label{eq:discbound}
\norm{\nabla u^h_k}_{\Lrm^6}&\leq C\, \norm{-\mathcal{L}_k^hu_k^h}_{\Lrm^2}&&\qquad\text{ if }d=3,\text{ and }\\
\norm{\nabla u^h_k}_{\Lrm^p}&\leq C\, \norm{-\mathcal{L}_k^hu_k^h}_{\Lrm^2}&&\qquad\text{ if }d=2\text{ for all }p\in (1,\infty).
\end{aligned}
\end{alignat}
The inequality~\eqref{eq:laplacelp} then implies that, for all $k\in \{0,\ldots,N\}$,
\begin{alignat}{2}
\norm{\nabla u^h_k}_{\Lrm^6}&\leq C\,\norm{f}_{\Lrm^\infty(\Lrm^2)}^{\max\{1,(q-1)/2\}}+C&&\qquad\text{ if }d=3,\text{ and } \label{eq:est_f1}
\\
\norm{\nabla u^h_k}_{\Lrm^p}&\leq C\,\norm{f}_{\Lrm^\infty(\Lrm^2)}^{\max\{1,(q-1)/2\}}+C&&\qquad\text{ if }d=2\text{ for all }p\in (1,\infty). \label{eq:est_f2}
\end{alignat}
Moreover, taking $q=4$, by Morrey's embedding theorem,
\begin{align}
\label{eq:space-d}
[u^h_k]_{\Crm^{0,\beta}}\leq  C\norm{f}_{\Lrm^\infty(\Lrm^2)}^{3/2}+C,
\end{align}
for some $\beta \in (0,1)$ and uniformly in $h,k$.

\subsection{Estimates for the time derivatives}

We test the $k$th inequality of~\eqref{eq:ELweak} with $\xi=0$ and the $(k-1)$st inequality with $u^h_k - u^h_{k-2}$, and divide by $\tau$ to find that
\begin{align*}
 &\int_\Omega R_1(\delta^{\tau,h}_k)+R_1(\delta_{k-1}^{\tau,h})+\nabla u^h_k : \Abb_{k} :  \nabla\delta^{\tau,h}_k-\nabla u^h_{k-1} : \Abb_{k-1} :  \nabla\delta^{\tau,h}_k
 \\
 &\qquad + \bigl[\DD W_0(u^h_k)-\DD W_0(u^h_{k-1})\bigr]\cdot\delta^{\tau,h}_k  -(f_k-f_{k-1})\cdot \delta^{\tau,h}_k  - R_1(\delta^{\tau,h}_k+\delta^{\tau,h}_{k-1}) \dd x \\
&\qquad \leq 0.
\end{align*}
Consequently,
\begin{align*}
 &\int_\Omega R_1(\delta^{\tau,h}_k)+R_1(\delta_{k-1}^{\tau,h})+(\nabla u^h_k-\nabla u^h_{k-1}) : \Abb_{k} :  \nabla\delta^{\tau,h}_k \\
&\qquad +[\DD W_0(u^h_k)-\DD W_0(u^h_{k-1})]\cdot\delta^{\tau,h}_k
-(f_k-f_{k-1})\cdot \delta^{\tau,h}_k  \\
&\qquad - R_1(\delta^{\tau,h}_k+\delta^{\tau,h}_{k-1}) - \nabla u^h_{k-1} : (\Abb_{k}-\Abb_{k-1}) :  \nabla\delta^{\tau,h}_k \dd x \leq 0.
\end{align*}
Dividing by $\tau > 0$ and using the subadditivity of $R_1$ together with the ellipticity gives
\begin{align*}
  \kappa \int_\Omega \abs{\nabla\delta^{\tau,h}_k}^2 \dd x
  &\leq \int_\Omega \frac{\abs{f_k-f_{k-1}}}{\tau} \, \abs{\delta^{\tau,h}_k}\dd x+\int_\Omega \frac{\abs{\Abb_{k}-\Abb_{k-1}}}{\tau}\, \abs{\nabla u^h_{k-1}}\, \abs{\nabla\delta^{\tau,h}_k}\dd x
  \\
  &\qquad- \int_\Omega \frac{1}{\tau^2} \bigl[ \DD W_0(u^h_k) - \DD W_0(u^h_{k-1}) \bigr]  \cdot (u^h_k - u^h_{k-1}) \dd x .
\end{align*}
The first term on the right-hand side can be bounded as
\begin{align}
\label{eq:fdot}
\begin{aligned}
  \int_\Omega \frac{\abs{f_k-f_{k-1}}}{\tau} \, \abs{\delta^{\tau,h}_k} \dd x
  &\leq \biggl( \dashint_{t_{k-1}}^{t_k}\norm{\dot{f}}_{\Lrm^2} \dd t \biggr) \norm{\delta^{\tau,h}_k}_{\Lrm^2} \\
  &\leq  C\biggl( \dashint_{t_{k-1}}^{t_k}\norm{\dot{f}}_{\Lrm^2} \dd t \biggr) \norm{\nabla \delta^{\tau,h}_k}_{\Lrm^2}.
  \end{aligned}
\end{align}
For the second term we deduce from~\eqref{eq:Lcoeff} and~\eqref{eq:h-coerc-d} that
\begin{align*}
\int_\Omega \frac{\abs{\Abb_{k}-\Abb_{k-1}}}{\tau}\, \abs{\nabla u^h_{k-1}} \, \abs{\nabla \delta^{\tau,h}_k}\dd x
&\leq C \norm{\Abb}_{\Crm^{0,1}([0,T] \times \Omega)} \, \norm{\nabla u^h_{k-1}}_{\Lrm^2}\,\norm{\nabla \delta^{\tau,h}_{k}}_{\Lrm^2} \\
&\leq C (1 + \norm{f_{k-1}}_{\Lrm^2})\,\norm{\nabla \delta^{\tau,h}_{k}}_{\Lrm^2}.
\end{align*}
For the third term we use~\eqref{eq:Wmon} and the Poincar\'{e} inequality to get
\begin{align*}
  - \frac{1}{\tau^2} \int_\Omega \bigl[ \DD W_0(u^h_k) - \DD W_0(u^h_{k-1}) \bigr]  \cdot (u^h_k - u^h_{k-1}) \dd x
  &\leq \mu \int_\Omega \abs{\delta^{\tau,h}_k}^2 \dd x \\
  &\leq \mu C_P(\Omega)^2 \int_\Omega \abs{\nabla \delta^{\tau,h}_k}^2 \dd x,
\end{align*}
where we recall that by $C_P(\Omega) > 0$ we denote the Poincar\'{e} constant of $\Omega$.
Hence, combining, we get
\begin{equation} \label{eq:apri_time-d}
 \norm{\nabla \delta^{\tau,h}_k}_{\Lrm^2} \leq \frac{C}{\kappa-\mu C_P(\Omega)^2} \biggl[ 1 + \dashint_{t_{k-1}}^{t_k}\norm{\dot{f}}_{\Lrm^2} \dd t +\norm{f_{k-1}}_{\Lrm^2} \biggr].
\end{equation}
By~\eqref{eq:muCP}, the constant appearing in front of the square bracket on the right-hand side is positive.

\subsection{Proof of Theorem~\ref{thm:approx}}
We define the approximation sequence
\[
  u^h_\tau(t):=\frac{t-t_{k-1}}{\tau}\,u^h_k+\frac{t_k-t}{\tau}\,u^h_{k-1}  \qquad\text{for $t \in (t_{k-1},t_k]$,\; $k = 1,\ldots,N$.}
\]
Integrating~\eqref{eq:apri_time-d} in time, we find (also using Jensen's inequality)
\begin{align*}
\norm{\nabla \dot{u}^h_\tau}_{\Lrm^a(\Lrm^2)}^a&=
\tau\sum_{k=0}^N\norm{\nabla \delta^{\tau,h}_k}_{\Lrm^2}^a
\leq C \tau\sum_{k=0}^N\biggl[ 1 +  \dashint_{t_{k-1}}^{t_k}\norm{\dot{f}}_{\Lrm^2}^a \dd t +\norm{f}_{\Lrm^\infty(\Lrm^2)}^a \biggr]
\\
&\leq C(1+\norm{\dot{f}}_{\Lrm^a(\Lrm^2)}^a+\norm{{f}}_{\Lrm^\infty(\Lrm^2)}^a).
\end{align*}
Now, by~\eqref{eq:space-d} and since $\dot{f}\in \Lrm^a(\Lrm^2)$ for some $a > 1$, we have that
\[
\sup_{t\in [0,T]}[u^h_\tau(t)]_{\Crm^{0,\beta}}+\norm{\nabla \dot{u}^h_{\tau}}_{\Lrm^a(\Lrm^2)}\leq C,
\]
uniformly in $h,\tau$. Hence, Remark~\ref{rem:H} yields that
\begin{equation} \label{eq:uh_apriori}
[u^h_\tau]_{\Crm^{0,\gamma}([0,T]\times \overline{\Omega})}+ \norm{\nabla \dot{u}^h_{\tau}}_{\Lrm^a(\Lrm^2)}\leq C,
\end{equation}
for some $\gamma \in (0,1)$.

Next we define $\tau_N=\frac{T}{N}$ and couple it with a sequence of space discretizations, with $h_N\to 0$,
as $N \rightarrow \infty$. Here, the $h_N$ is associated with the approximation space $\Xrm^{h_N}$.
Let $u_N:=u^{h_N}_{\tau_N}$ be the sequence of discrete solutions to \eqref{eq:discrete1-d} with respect to
$\tau_N$ and $\Xrm^{h_N}$. The above estimates imply that there is a subsequence of the sequence of
approximate solutions converging uniformly in $\Crm^{0,\gamma_0}([0,T]\times \Omega;\R^m)$  by the Arzel\`a--Ascoli theorem,
for some $\gamma_0 \in (0,\gamma)$, and weakly in $\Wrm^{1,a}(0,T;\Wrm^{1,2}_0(\Omega;\R^m))$.
The proof that this subsequence, in fact, converges to the (unique) solution of the problem proceeds as in the proof
of the existence of solutions in Section~\ref{ssc:proof}, except for one difference: since one cannot use a fixed
function $\xi\in \Lrm^1(0,T;\Wrm^{1,2}_0(\Omega;\R^m))$ as a test function on the discrete level, one has to replace it
with a sequence of admissible functions $(\xi^N)_{N\in \Nbb}$ such that
$\xi^N\to \xi$, in $\Lrm^1(0,T;\Wrm^{1,2}_0(\Omega;\R^m))$.
For example, one can take $\xi^N$ to be continuous piecewise affine function with respect to $t$ on the partition $\{t_0^N,\dots,t_N^N\}$
of $[0,T]$ such that $\xi^N(t) \in \Xrm^{h_N}$ for all $t \in [0,T]$, and then appeal to the fact that
$\bigcup\limits_{N\in \Nbb} \Xrm^{h_N}$ is dense in $\Wrm^{1,2}_0(\Omega;\R^m)$,
and the set of all such piecewise affine functions defined on $[0,T]$ is dense in $\Lrm^1(0,T)$.

\section{Convergence with a rate} \label{sc:rate}

The aim of this section is to derive a bound on the error between the analytical solution to the model and its fully discrete approximation constructed in the previous section, under the hypotheses adopted on the data. Our main result in this respect is encapsulated in the following theorem.

\begin{theorem} \label{thm:rate}
Under the assumptions made on $\Lcal,f,u_0$ and $W_0$, the sequence of approximations generated by the numerical method formulated in the previous section
converges to the solution $u$ with a rate, in the sense that there exists a positive constant $C$, independent of $h$ and $\tau$, such that
\[
\int_0^T\norm{u^h_{\tau}-u}_{\Wrm^{1,2}}^2 \dd t\leq C\, (h +\tau^{\min\{1,a-1\}}),
\]
where $a\in (1,\infty]$ is the assumed integrability exponent of $\dot{f}$ in time.
\end{theorem}

We embark on the proof of Theorem \ref{thm:rate} by developing, in the next two subsections, some preliminary results that are required in the proof of the final error bound.

\subsection{Key estimates that are used in the proof of Theorem~\ref{thm:rate}}
We will collect here the bounds established in earlier sections that are critical for the proof of Theorem~\ref{thm:rate}.

Thanks to the a-priori estimates leading to~\eqref{eq:uh_apriori}, we have that
\[
\norm{u^h_{\tau}}_{\Lrm^\infty((0,T)\times\Omega)}
+\norm{\nabla {\dot{u}^h_{\tau}}}_{\Lrm^a(\Lrm^2)}\leq C,
\]
uniformly in $h$ and $\tau$; here we have also used that $q \leq 4$ if $d = 3$.
Moreover, thanks to Theorem~\ref{thm:exist} we have that
\[
\norm{u}_{\Lrm^\infty((0,T)\times\Omega)}
+\norm{\nabla {\dot{u}}}_{\Lrm^a(\Lrm^2)}\leq C.
\]
These bounds will be used repeatedly throughout the proof without mentioning them explicitly.
Without loss of generality we will assume that $h,\tau\in (0,1)$.
Also, as is clear from the statement of Theorem~\ref{thm:rate}, there is no improvement
in the rate of convergence with respect to $\tau$ for $a>2$;
we shall therefore assume in the proof of Theorem~\ref{thm:rate} that $a\in (1,2]$.

\subsection{A C\'ea type inequality}
We recall the definition
\[
  u^h_{\tau}(t):=\frac{t-t_{k-1}}{\tau}\,u^h_k+\frac{t_k-t}{\tau}\,u^h_{k-1}  \qquad\text{for $t \in (t_{k-1},t_k]$,\; $k = 1,\ldots,N$.}
\]
In the following we will fix a $t\in (t_{k-1},t_{k}]$, such that $\norm{\nabla \dot{u}(t)}_{\Lrm^2}<\infty$, which holds for almost every $t\in (t_{k-1},t_{k}]$. Henceforth, all time-dependent functions considered on $(t_{k-1},t_k]$ will be evaluated at that point $t$.

We start with testing the continuous system. We use the test function $\delta^{\tau,h}_k=\dot{u}^h_{\tau}$, which is admissible by~\eqref{eq:apri_time-d}.
Then we find for the chosen $t\in (t_{k-1},t_{k}]$, that
\begin{align}
\label{est:1}
\begin{aligned}
&\int_\Omega R_1(\dot{u}) + \nabla u : \Abb : \nabla(\dot{u}-\dot{u}^h_{\tau})+\bigl[\DD W_0(u) - f \bigr]\cdot(\dot{u}-\dot{u}^h_{\tau}) \dd x
\leq \int_\Omega R_1(\dot{u}^h_{\tau}) \dd x.
\end{aligned}
\end{align}
We now wish to imitate this inequality for the discrete problem.
We introduce for each time step $k \in \{1,\ldots,N\}$ the Ritz projector
\[
\Pi^h_k: \Wrm^{1,2}_0(\Omega;\R^m)\to \Xrm^h\text{ defined as }\skp{\Abb_{k}\nabla (g-\Pi^h_k g)}{\nabla \phi^h}=0\text{ for all }\phi^h\in \Xrm^h.
\]
For any $g\in \Wrm^{1,2}_0(\Omega;\R^m)$ we have that
\begin{align*}
\norm{\nabla (g-\Pi^h_k g)}_{\Lrm^2}^2&\leq C\skp{\Abb_k\nabla (g-\Pi^h_k g)}{\nabla (g-\Pi^h_k g)}\\
&=C\skp{\Abb_k\nabla (g-\Pi^h_k g)}{\nabla (g-\Pbb^hg)}
\\
&\leq C \norm{\nabla (g-\Pi^h_k g)}_{\Lrm^2}\norm{\nabla (g-\Pbb^h g)}_{\Lrm^2},
\end{align*}
where we have used the definition of $\Pi^h_k$ in the transition to the second line.
Hence, using the properties of the projector $\Pbb^h$, it follows that
\begin{alignat}{2}
\label{eq:proj1-k}
\norm{\nabla (g-\Pi^h_k g)}_{\Lrm^2}
&\leq C \norm{g}_{\Wrm^{1,2}}&&\qquad\text{ for all $g\in \Wrm^{1,2}_0(\Omega;\R^m)$},
\\
\label{eq:proj2-k}
\norm{\nabla (g-\Pi^h_k g)}_{\Lrm^2}
&\leq Ch\norm{g}_{\Wrm^{2,2}}&&\qquad\text{ for all $g\in (\Wrm^{2,2}\cap \Wrm^{1,2}_0)(\Omega;\R^m)$},
\end{alignat}
with a constant $C > 0$ that depends only on the ellipticity of $\Abb$ and the properties of $\Pbb$, but is independent
of the discretization parameters $h$ and $\tau$. We also require a bound in the $\Lrm^2$ norm on
the difference between a function $g \in (\Wrm^{2,2} \cap \Wrm^{1,2}_0)(\Omega;\R^m)$ and its Ritz projection $\Pi^h_k g$.
The derivation of this proceeds by an Aubin--Nitsche type duality argument: for
$g \in (\Wrm^{2,2} \cap \Wrm^{1,2}_0)(\Omega;\R^m)$ we let $v_g\in \Wrm^{1,2}_0(\Omega;\R^m)$ be the unique solution of the problem
 \[
 \skp{ \nabla v_g}{\Abb_k\nabla \phi}=\skp{\Abb_k \nabla v_g}{\nabla \phi}= \skp{g-\Pi^h_k g}{\phi}\quad \text{ for all }\phi\in \Wrm^{1,2}_0(\Omega;\R^m).
 \]
Hence, by the regularity theory of uniformly elliptic systems,~\eqref{eq:proj1-k} and~\eqref{eq:proj2-k} yield that
 \begin{align*}
 &\norm{g-\Pi^h_k g}_{\Lrm^2}^2=\skp{\nabla v_g}{\Abb_k \nabla(g-\Pi^h_k g)}\\
 &\qquad=\skp{\nabla( v_g-\Pi^h_k v_g)}{\Abb_k \nabla(g-\Pi^h_k g)}
 \\
&\qquad \leq C \norm{\nabla (v_g-\Pi^h_k v_g)}_{\Lrm^2}\cdot \norm{\nabla(g-\Pi^h_k g )}_{\Lrm^2}\\
&\qquad \leq Ch \norm{v_g}_{\Wrm^{2,2}} \cdot h^s \|g\|_{\Wrm^{1+s,2}} \\
&\qquad \leq Ch^{1+s} \norm{g-\Pi^h_k g}_{\Lrm^{2}}\|g\|_{\Wrm^{1+s,2}}, \qquad s=0,1.
 \end{align*}
Therefore,
\begin{align}
 \label{eq:ritz0}
 \begin{aligned}
\norm{g-\Pi^h_k g}_{\Lrm^2} \leq C h^{1+s} \|g\|_{\Wrm^{1+s,2}}, \qquad s=0,1.
\end{aligned}
\end{align}

Now, for the fixed $t\in (t_{k-1},t_{k}]$, we find
\begin{align}
\label{eq:preest1}
\begin{aligned}
&\int_\Omega R_1(\dot{u}^h_{\tau}) + \nabla u^h_{\tau} : \Abb_{k} : \nabla(\dot{u}^h_{\tau}-\dot{u})\dd x \\
&\qquad= \int_\Omega R_1(\dot{u}^h_{\tau}) + \nabla u^h_{\tau} : \Abb_{k} : \nabla(\dot{u}^h_{\tau}-\Pi^h_k\dot{u}) \dd x
\\
&\qquad = \int_\Omega R_1(\dot{u}^h_{\tau}) +\nabla u_k^h : \Abb_{k} : \nabla(\dot{u}^h_{\tau}-\Pi^h_k\dot{u})\dd x
\\
&\qquad\qquad +\int_\Omega \frac{t_k-t}{\tau}\nabla(u_{k-1}^h- u_k^h) : \Abb_{k} : \nabla(\dot{u}^h_{\tau}-\Pi^h_k\dot{u})\dd x,
\end{aligned}
\end{align}
using the fact that
\[
  u^h_{\tau}(t):=u^h_{k}+\frac{t_{k}-t}{\tau}\,(u^h_{k-1}-u^h_{k})  \qquad\text{for $t \in (t_{k-1},t_k]$,\; $k = 1,\ldots,N$.}
\]
We now proceed to estimate using~\eqref{eq:proj1-k}, recalling that $a\in (1,2]$ and using the fact
that $\nabla u_k^h$ is bounded in $\Lrm^2$, uniformly in $k,h$, thanks
to~\eqref{eq:est_f1},~\eqref{eq:est_f2}. Hence,
\begin{align}
&\frac{t_k-t}{\tau}\int_\Omega\nabla(u_{k-1}^h- u_k^h) : \Abb_{k} : \nabla(\dot{u}^h_{\tau}-\Pi^h_k\dot{u})\dd x
\notag\\
&\qquad \leq
C \norm{\nabla(u_{k-1}^h- u_k^h)}_{\Lrm^2}\norm{\nabla(\dot{u}^h_{\tau}-\Pi^h_k\dot{u})}_{\Lrm^2}
\notag\\
&\qquad \leq C \bigl( \norm{\nabla(u_{k-1}^h- u_k^h)}_{\Lrm^2}\norm{\nabla\dot{u}^h_{\tau}}_{\Lrm^2}+\norm{\nabla(u_{k-1}^h- u_k^h)}_{\Lrm^2}\norm{\nabla (\Pi^h_k\dot{u})}_{\Lrm^2} \bigr)
\notag\\
&\qquad \leq C\tau^{a-1}\norm{\nabla(u_{k-1}^h- u_k^h)}_{\Lrm^2}^{2-a}\norm{\nabla\dot{u}^h_{\tau}}_{\Lrm^2}^a
+
C\tau^{a-1}\norm{\nabla(u_{k-1}^h- u_k^h)}_{\Lrm^2}^{2-a}\norm{\nabla\dot{u}}_{\Lrm^2}^a
\notag\\
&\qquad \leq  C\tau^{a-1}(\norm{\nabla\dot{u}^h_{\tau}}_{\Lrm^2}^a+\norm{\nabla\dot{u}}_{\Lrm^2}^a).
\label{eq:preest2}
\end{align}
Here in the transition from the third line to the fourth line we have noted that
$\dot{u}^h_\tau(t)=\frac{u^h_k-u^h_{k-1}}{\tau}$ at the point $t \in (t_{k-1},t_k]$,
and used the stability of the Ritz projector
$\Pi^h_k$ in the $\Wrm^{1,2}(\Omega;\R^m)$-seminorm together with Young's inequality,
to deduce the following bound on the second term in the third line:
\begin{align*}
\norm{\nabla(u_{k-1}^h- u_k^h)}_{\Lrm^2}\norm{\nabla (\Pi^h_k\dot{u})}_{\Lrm^2} &\leq
C\,\norm{\nabla(u_{k-1}^h- u_k^h)}_{\Lrm^2}\norm{\nabla \dot{u}}_{\Lrm^2}\\
&= C \tau^{a-1}\norm{\nabla(u_{k-1}^h- u_k^h)}_{\Lrm^2}^{2-a}\norm{\nabla \dot{u}^h_\tau}_{\Lrm^2}^{a-1}\norm{\nabla \dot{u}}_{\Lrm^2}
\\
& \leq C\tau^{a-1}\norm{\nabla(u_{k-1}^h- u_k^h)}_{\Lrm^2}^{2-a} \bigl(\norm{\nabla \dot{u}^h_\tau}_{\Lrm^2}^{a}+\norm{\nabla \dot{u}}_{\Lrm^2}^a\bigr),
\end{align*}
and the first part of this expression can be combined with the first part in the fourth line of~\eqref{eq:preest2}.

Next we apply~\eqref{eq:discrete1-d} with the test function $\xi^h=\Pi^h_k\dot{u}$.
This implies, noting that
$\dot{u}^h_\tau(t)=\frac{u^h_k-u^h_{k-1}}{\tau}=\delta^{\tau,h}_k$ at the point $t \in (t_{k-1},t_k]$, and using first \eqref{eq:preest1} and then \eqref{eq:preest2} that
\begin{align}
&\hspace{-1cm}\int_\Omega R_1(\dot{u}^h_{\tau}) + \nabla u^h_{\tau} : \Abb_{k} : \nabla(\dot{u}^h_{\tau}-\dot{u})\dd x
\notag\\
&\hspace{-3mm}\leq \int_\Omega R_1(\Pi^h_k \dot{u})+ \bigl[ -\DD W_0(u^h_k)
+ f_k \bigr]\cdot(\dot{u}^h_{\tau}-\Pi^h_k\dot{u}) \notag\\
&\qquad\qquad +\frac{t_k-t}{\tau}\nabla(u_{k-1}^h- u_k^h) : \Abb_{k} : \nabla(\dot{u}^h_{\tau}-\Pi^h_k\dot{u})\dd x
\notag\\
&\hspace{-3mm}\leq \int_\Omega R_1(\Pi^h_k \dot{u})+\bigl[ -\DD W_0(u^h_k) + f_k \bigr]\cdot(\dot{u}^h_{\tau}-\Pi^h_k\dot{u}) \dd x \notag\\
&\qquad + C\tau^{a-1}(\norm{\nabla\dot{u}^h_{\tau}}_{\Lrm^2}^a+\norm{\nabla\dot{u}}_{\Lrm^2}^a).  \label{est:2}
\end{align}
As $\nabla u$ in bounded in $\Lrm^\infty(0,T;\Lrm^2(\Omega;\R^{d\times m}))$, $\nabla u^h_{\tau}$ is uniformly bounded (with respect to $h$ and $\tau$) in $\Lrm^\infty(0,T;\Lrm^2(\Omega;\R^{d\times m}))$, and since $\Abb$ is assumed to be Lipschitz continuous, we have that
\begin{align*}
(I)&:=\int_\Omega \nabla(u- u^h_{\tau}) : \Abb : \nabla(\dot{u}-\dot{u}^h_{\tau})\dd x
\\
&=\int_\Omega \nabla u:\Abb: \nabla(\dot{u}-\dot{u}^h_{\tau}) + \nabla u^h_{\tau} : \Abb_{k} : \nabla(\dot{u}^h_{\tau}-\dot{u}) \\
&\qquad\qquad + \nabla u^h_{\tau} :(\Abb-\Abb_{k}) : \nabla(\dot{u}^h_{\tau}-\dot{u})\dd x
\\&\leq\int_\Omega \nabla u:\Abb: \nabla(\dot{u}-\dot{u}^h_{\tau}) + \nabla u^h_{\tau} : \Abb_{k} : \nabla(\dot{u}^h_{\tau}-\dot{u})\dd x+ C\tau(\norm{\nabla \dot{u}^h_\tau}_{\Lrm^2}+\norm{\nabla\dot{u}}_{\Lrm^2}).
\end{align*}
This allows us to use the estimate derived by summing up~\eqref{est:1} and~\eqref{est:2}, which implies that
\begin{align}\label{eq:1234}
(I)
&\leq \int_\Omega R_1(\Pi^h_k \dot{u})-R_1(\dot{u})\dd x
 + \int_\Omega \DD W_0 (u^h_k)\cdot (\Pi^h_k \dot{u}-\dot{u}^h_{\tau})+\DD W_0(u)\cdot(\dot{u}^h_{\tau}-\dot{u})\dd x
 \notag\\
 &\qquad
  + \int_\Omega f_k\cdot (\dot{u}^h_{\tau}-\Pi^h_k \dot{u})+f \cdot (\dot{u}-\dot{u}^h_{\tau})\dd x +C\tau^{a-1}(\norm{\nabla\dot{u}^h_{\tau}}_{\Lrm^2}^a+\norm{\nabla\dot{u}}_{\Lrm^2}^a)
  \notag\\
  &\qquad
 =: (II)+(III)+(IV)+C\tau^{a-1}(\norm{\nabla\dot{u}^h_{\tau}}_{\Lrm^2}^a+\norm{\nabla\dot{u}}_{\Lrm^2}^a).
 \end{align}

\subsection{Imitation of the uniqueness estimate for the error}
Using the 1-homogeneity and subadditivity of $R_1$, we obtain
 \begin{align}
 (II)\leq C \int_{\Omega}\abs{\dot{u}-\Pi^h_k \dot{u}}\dd x \leq C \norm{\dot{u}-\Pi^h_k\dot{u}}_{\Lrm^2}.
 \end{align}
Consequently, by taking $s=0$, $g=\dot{u}$ in~\eqref{eq:ritz0}, and using Poincar\'e's inequality,
\begin{align}
\label{eq:II-d}
(II)\leq C \norm{\dot{u}-\Pi^h_k \dot{u}}_{\Lrm^2}\leq  Ch \norm{\nabla \dot{u}}_{\Lrm^2}.
\end{align}

Next, using the fact that $\abs{\DD W_0(u^h_{\tau})}$ is uniformly bounded and noting~\eqref{eq:II-d}, we find that
 \begin{align*}
 (III)&=\int_\Omega \bigl[\DD W_0 (u^h_k)-\DD W_0(u^h_{\tau})\bigr]\cdot (\dot{u}-\dot{u}^h_{\tau})+
 \bigl[\DD W_0(u)-\DD W_0(u^h_{\tau})\bigr]\cdot (\dot{u}^h_{\tau}-\dot{u})
 \\&\qquad+\DD W_0(u^h_k)\cdot(\Pi^h_k \dot{u}-\dot{u})\dd x
\\
 &\leq \int_\Omega \bigl[\DD W_0 (u^h_k)-\DD W_0(u^h_{\tau})\bigr]\cdot (\dot{u}-\dot{u}^h_{\tau})+
 \!\bigl[\DD W_0(u)-\DD W_0(u^h_{\tau})\bigr]\cdot (\dot{u}^h_{\tau}-\dot{u})\dd x \\
&\qquad +Ch \norm{\nabla \dot{u}}_{\Lrm^2}.
 \end{align*}
Moreover, using that  $\DD^2 W_0(\theta u^h_k+(1-\theta) u^h_{\tau})$ is uniformly bounded with respect to $h$ and $\tau$ thanks to the uniform bounds on $u^h_{\tau}$ and the assumptions on $W_0$,
noting that $\abs{u^h_{\tau}(t)-u^h_k(t)}\leq \tau \abs{\dot{u}^h_\tau}$ for all $t\in (t_{k-1}, t_{k-1}]$, and that, by hypothesis $a \in (1,2]$, we find that, pointwise,
\begin{align*}
\abs{\DD W_0 (u^h_k)-\DD W_0(u^h_{\tau})}&\leq\absBB{\int_0^1 \DD^2 W_0(\theta u^h_k+(1-\theta) u^h_{\tau})\dd \theta} \, \abs{u^h_k-u^h_{\tau}}
\\
&\leq C\abs{u^h_k-u^h_{\tau}}\\
&\leq C \tau^{a-1}\norm{u_h^\tau}_{\Lrm^\infty}^{2-a}\abs{\dot{u}^h_{\tau}}^{a-1}\\
&\leq C \tau^{a-1}\abs{\dot{u}^h_{\tau}}^{a-1}.
\end{align*}
 Therefore, combining the above estimates with Young's inequality (using $\frac{1}{a}+\frac{a-1}a=1$) and the Sobolev--Poincar\'e inequality implies that
 \begin{align*}
 (III)
 &\leq \int_\Omega \bigl[\DD W_0 (u)-\DD W_0(u^h_{\tau})\bigr]\cdot (\dot{u}^h_{\tau}-\dot{u})\dd x
 + C\tau^{a-1}\int_{\Omega}\abs{\dot{u}^h_\tau}^{a-1}\abs{\dot{u}-\dot{u}^h_\tau}\dd x \\
&\qquad +Ch \norm{\nabla \dot{u}}_{\Lrm^2}
 \\
&\leq \int_\Omega \bigl[\DD W_0 (u)-\DD W_0(u^h_{\tau})\bigr]\cdot (\dot{u}^h_{\tau}-\dot{u})\dd x+ C\tau^{a-1}(\norm{\nabla \dot{u}^h_\tau}_{\Lrm^2}^a+\norm{\nabla \dot{u}}_{\Lrm^2}^a)\\
&\qquad +Ch\norm{\nabla \dot{u}}_{\Lrm^2}.
 \end{align*}
 Using the symmetry of $\DD^2 W_0$ yields that
 \begin{align*}
 &\bigl[\DD W_0 (u)-\DD W_0(u^h_{\tau})\bigr]\cdot (\dot{u}^h_{\tau}-\dot{u})\\
&\qquad = \int_0^1\DD^2 W_0(\theta u+(1-\theta)u^h_\tau )\dd \theta\,
[u-u^h_\tau ,\dot{u}^h_\tau-\dot{u}]
\\
&\qquad = -\int_0^1\DD^2 W_0(\theta u+(1-\theta)u^h_\tau )\dd \theta\,
[u-u^h_\tau ,\dot{u}-\dot{u}^h_\tau]
\\
&\qquad = -\frac{\di}{\di t} \bigg(\frac12\int_0^1\DD^2 W_0(\theta u+(1-\theta)u^h_\tau )\dd \theta\,
[u-u^h_\tau ,u-u^h_\tau]\bigg)
\\
&\qquad \qquad +\frac12\int_0^1\partial_t \Big(\DD^2 W_0(\theta u+(1-\theta)u^h_\tau ))\dd \theta\Big)\, [u-u^h_\tau ,u-u^h_\tau].
 \end{align*}
 We estimate further using~\eqref{eq:DDDW} and the a-priori bounds in the space-time $\Lrm^\infty$ norm on $u^h_{\tau}$ and $u$, H\"older's inequality and the Sobolev--Poincar\'{e} inequality, to get
 \begin{align*}
&\int_\Omega  \int_0^1\partial_t \Big(\Big(\DD^2 W_0(\theta u+(1-\theta)u^h_\tau ))\Big)\dd \theta \,[u-u^h_\tau ,u-u^h_\tau]\dd x
\\
 &\qquad \leq \int_0^1 \abs{\DD^3(\theta u+(1-\theta)u^h_\tau )) \, \abs{\theta\dot{u}+(1-\theta)\dot{u}^h_\tau}\,\abs{u-u^h_\tau}^2} \dd \theta
 \\
&\qquad \leq C(\norm{\dot{u}}_{\Lrm^2}+\norm{\dot{u}^h_\tau}_{\Lrm^2})\,\norm{\nabla (u-u^h_\tau)}_{\Lrm^2}^2.
 \end{align*}
 Finally we get that
 \begin{align*}
 (III)
 &\leq -\frac{\di}{\di t} \biggl(\int_\Omega \frac12\int_0^1\DD^2 W_0(\theta u+(1-\theta)u^h_\tau )\dd \theta\,
[u-u^h_\tau ,u-u^h_\tau] \dd x\biggr)
\\
&\qquad+ C(\norm{\dot{u}}_{\Lrm^2}+\norm{\dot{u}^h_\tau}_{\Lrm^2})\,\norm{\nabla (u-u^h_\tau)}_{\Lrm^2}^2 +
C\tau^{a-1}(\norm{\nabla \dot{u}^h_\tau}_{\Lrm^2}^a+\norm{\nabla \dot{u}}_{\Lrm^2}^a) \\
&\qquad +Ch\norm{\nabla \dot{u}}_{\Lrm^2}.
 \end{align*}
The last term is estimated using that, for $a \in (1,2]$ and $t\in(t_{k-1},t_k]$,
\[
\abs{f(t_k)-f(t)}\leq 2^{2-a} \norm{f}_{\Lrm^\infty}^{2-a}\absBB{\int_{t}^{t_k}\dot{f}\dd s}^{a-1}\leq C\tau^{a-1}\bigg(\dashint_{t_{k-1}}^{t_k}\abs{\dot{f}}\dd s\bigg)^{a-1}.
\]
Hence, using Young's inequality, the Sobolev--Poincar\'{e} inequality, \eqref{eq:proj1-k}, Jensen's inequality, Fubini's theorem and \eqref{eq:II-d}, we find
\begin{align*}
 (IV)&= \int_\Omega (f_k-f)\cdot (\dot{u}^h_{\tau}-\Pi^h_k\dot{u})+f \cdot (\dot{u}-\Pi^h_k\dot{u})\dd x
 \\
 &\leq C\tau^{a-1}\int_\Omega \bigg(\dashint_{t_{k-1}}^{t_k}\abs{\dot{f}}\dd s\bigg)^{a-1}(\abs{\dot{u}^h_{\tau}}+\abs{\Pi^h_k\dot{u}})\dd x+Ch\norm{f}_{\Lrm^\infty} \norm{\nabla \dot u}_{\Lrm^2}
 \\
 &\leq C\tau^{a-1}\bigg(\int_\Omega \bigg(\dashint_{t_{k-1}}^{t_k}\abs{\dot{f}}\dd s\bigg)^{a}\dd x+\norm{\nabla\dot{u}}_{\Lrm^2}^a+\norm{\nabla \dot{u}_\tau^h}_{\Lrm^2}^a\bigg)+Ch\norm{\nabla \dot u}_{\Lrm^2}
 \\
 &\leq C\tau^{a-1}\bigg(\dashint_{t_{k-1}}^{t_k}\norm{\dot{f}}_{\Lrm^a}^a\dd s+\norm{\nabla\dot{u}}_{\Lrm^2}^a+\norm{\nabla \dot{u}_\tau^h}_{\Lrm^2}^a\bigg)+Ch\norm{\nabla \dot u}_{\Lrm^2}.
\end{align*}
We now wish to replicate the Gronwall argument from the uniqueness section. We therefore define the error function
\[
e(t):=\norm{\nabla(u-u^h_{\tau})}_{\Lrm^2}^2.
\]
Using the uniform ellipticity assumption on $\Abb$, we find that
 \begin{align*}
(I')&:=\frac{\di}{\di t}\int_\Omega \nabla(u- u^h_{\tau}) : \frac{\Abb}{2} : \nabla(u-u^h_{\tau})\dd x\\
&=\int_\Omega \nabla(u- u^h_{\tau}) : \Abb : \nabla(\dot{u}-\dot{u}^h_{\tau})\dd x + \int_\Omega \nabla(u- u^h_{\tau}) : \frac{\partial_t{\Abb}}{2} : \nabla({u}-{u}^h_{\tau})\dd x
 \\
 &\leq (I) + C \int_\Omega\abs{\nabla (u-u^h_{\tau})}^2\dd x \\
 &\leq (I)+C e(t).
 \end{align*}

 Collecting the bounds on the terms $(II)$, $(III)$ and $(IV)$ and substituting them into~\eqref{eq:1234}, yields a bound on $(I)$, from which we then deduce that
 \begin{align*}
 &(I')+\frac{\di}{\di t} \bigg(\int_\Omega \frac12\int_0^1\DD^2 W_0(\theta u(t)+(1-\theta)u^h_\tau(t) )\dd \theta\,
[u-u^h_\tau ,u-u^h_\tau] \dd x\bigg)
\\
&\qquad \leq C(h+\tau^{a-1})\bigg(1+ \ONE_{(t_{k-1},t_k]}(t) \dashint_{t_{k-1}}^{t_k}\norm{\dot{f}}_{\Lrm^a}^a\dd s+\norm{\nabla \dot{u}^h_\tau(t)}_{\Lrm^2}^a+\norm{\nabla \dot{u}(t)}_{\Lrm^2}^a\bigg)
\\
&\qquad \qquad  + C(1+\norm{\dot{u}(t)}_{\Lrm^2}+\norm{\dot{u}^h_\tau(t)}_{\Lrm^2})
  \,e(t),\qquad t \in (t_{k-1},t_k],\quad k=1,\dots,N.
 \end{align*}
For any $s \in (0,T]$ there exists a $k \in \{1,\dots,N\}$ such that $s \in (t_{k-1},t_k]$. If $k=1$ we integrate
the above inequality over $(0,s]$; else we integrate it over each of the intervals $(0,t_1],\dots,(t_{k-1},s]$,
and sum up the resulting inequalities. Thus we deduce that
   \begin{align*}
&\int_{\Omega} \nabla(u(s)- u^h_{\tau}(s)) : \frac{\Abb(s)}{2} : \nabla({u}(s)-{u}^h_{\tau}(s))
\\
&\qquad \quad + \frac12\int_0^1\DD^2 W_0(\theta u(s)+(1-\theta)u^h_{\tau}(s))\dd \theta\,
[u(s)-u^h_{\tau}(s),u(s)-u^h_{\tau}(s)]\dd x
 \\
&\qquad \leq C \int_0^s(1+\norm{\dot{u}}_{\Lrm^2}+\norm{\dot{u}^h_\tau}_{\Lrm^2}) e(t)\dd t
\\
&\qquad \quad
\, + \,C(h+\tau^{a-1})(1+\norm{\dot{f}}_{\Lrm^a(\Lrm^a)}^a+\norm{\nabla \dot{u}^h_\tau}_{\Lrm^a(\Lrm^2)}^{a}+\norm{\nabla \dot{u}}_{\Lrm^{a}(\Lrm^2)}^{a})
\\
&\qquad \quad + \int_{\Omega} \nabla(u_0- u_0^h) : \frac{\Abb_0}{2} : \nabla(u_0-u_0^h)
\\
&\qquad \quad + \frac12\int_0^1\DD^2 W_0(\theta u_0+(1-\theta)u_0^h)\dd \theta\,
[u_0-u_0^h,u_0-u_0^h]\dd x.
 \end{align*}
We use~\eqref{eq:el-u0},~\eqref{eq:ELweak},~\eqref{eq:uh0} and~\eqref{eq:DWgrowth} to find that
\begin{align*}
(I'')&:=\int_{\Omega} \nabla(u_0- u_0^h) :\Abb_0 : \nabla(u_0-u_0^h)
\\
&\qquad + \int_0^1\DD^2 W_0(\theta u_0+(1-\theta)u_0^h)\dd \theta\,
[u_0-u_0^h,u_0-u_0^h]\dd x
\\
& = \int_{\Omega} \nabla u_0 :\Abb_0 : \nabla(u_0-u_0^h)+\DD W_0(u_0)\cdot [u_0-u_0^h]\dd x
\\
&\qquad + \int_{\Omega} \nabla u^h_0 :\Abb_0 : \nabla(u_0^h-u_0)+\DD W_0(u_0^h)\cdot [u_0^h-u_0]\dd x
\\
& \leq \int_{\Omega} R_1(u_0^h-u_0) + f(0)\cdot (u_0-u_0^h)\dd x
\\
&\qquad + \int_{\Omega} \nabla u^h_0 :\Abb_0 : \nabla(u_0^h-\Pi^h_0 u_0)+\DD W_0(u_0^h)\cdot [u_0^h-\Pi^h_0 u_0]\dd x
\\
&\qquad +\int_{\Omega} \nabla u^h_0 :\Abb_0 : \nabla(\Pi^h_0 u_0 -u_0) \dd x
+ \int_\Omega \DD W_0(u_0^h)\cdot [\Pi^h_0 u_0-u_0]\dd x
\\
& \leq \int_{\Omega} R_1(\Pi^h_0 u_0-u_0)\dd x+\big(\norm{f(0)}_{\Lrm^\infty} + C(1+\norm{u^h_0}_{\Lrm^\infty}^{q-1})\big)\int_{\Omega}\abs{\Pi^h_0 u_0-u_0}\dd x.
\end{align*}
Here in the transition from the right-hand side of the first
inequality on the right-hand side of the second
inequality we have used that the second integral on the right-hand side of
the first inequality is equal to zero since $u^h_0=\Pi^h_0u_0$, and
the third integral on the right-hand side of the first inequality
is also equal to zero as $u^h_0=\Pi^h_0u_0$ is the Ritz projection of $u_0$.

Next, using that $u_0,u_0^h\in (\Wrm^{1,2}_0\cap \Lrm^\infty)(\Omega;\R^m)$ by~\eqref{eq:space-d}, we find by~\eqref{eq:ritz0} that
\[
(I'')\leq C \int_{\Omega}\abs{\Pi^h_0(u_0)-u_0}\dd x\leq Ch\norm{\nabla u_0}_{\Lrm^2}.
\]
Hence, by using~\eqref{eq:gw3}, the bounds on $\nabla \dot{u},\nabla \dot{u}^h_\tau$, the assumptions on $\dot{f}$, and also the (uniform) positive-definiteness of $\Abb(s)$, we have that
\[
e(s)\leq C \int_0^s(1+\norm{\dot{u}}_{\Lrm^2}+\norm{\dot{u}^h_\tau}_{\Lrm^2})\, e(t)\dd t +  C(h+\tau^{a-1}).
\]
Since $1+\norm{\dot{u}}_{\Lrm^1(\Lrm^2)}+\norm{\dot{u}^h_\tau}_{\Lrm^1(\Lrm^2)}$ is bounded (uniformly in $h,\tau$) Gronwall's lemma implies the desired error bound. That completes the proof of Theorem \ref{thm:rate}.\hfill $\Box$

\appendix

\section*{Appendix}

Our aim in this appendix is to prove the inequalities~\eqref{eq:discbound}.  Let $z^h \in \Xrm^h$,
and define the discrete Green operator $\mathcal{G}^h_k\,: \Xrm^h \rightarrow \Xrm^h$ by
\[ \int_\Omega \nabla(\mathcal{G}^h_k z^h) : \mathbb{A}_k : \nabla \hat{\xi}^h \dd x = \int_\Omega z^h \cdot \hat{\xi}^h \dd x  \qquad \text{for all $\hat{\xi}^h \in \Xrm^h$,}\]
which is the finite element counterpart of the operator  $\mathcal{G}_k\,: \Lrm^2(\Omega;\R^m) \rightarrow \Wrm^{1,2}_0(\Omega;\R^m)$, defined, for any
$z \in \Lrm^2(\Omega;\R^m)$, by
\[ \int_\Omega \nabla(\mathcal{G}_k z) : \mathbb{A}_k : \nabla \phi \dd x = \int_\Omega z \cdot \phi \dd x  \qquad \text{for all $\phi  \in \Wrm^{1,2}_0(\Omega;\R^m)$.} \]
By elliptic regularity theory $\mathcal{G}_k z \in (\Wrm^{2,2} \cap \Wrm^{1,2}_0)(\Omega;\R^m)$ for all $k \in \{0,\dots,N\}$, and
\[ \|\mathcal{G}_k z\|_{\Wrm^{2,2}} \leq C \|z\|_{\Lrm^2},\]
where $C$ is a positive constant, independent of $k$.

It is easily seen that $\mathcal{G}^h_k \mathcal{L}^h_k z^h = z^h$ for all $z^h \in \Xrm^h$. Indeed, using the definitions of $\mathcal{G}^h_k$ and $\mathcal{L}^h_k$, we have that
\begin{eqnarray*}
\int_\Omega \nabla(\mathcal{G}^h_k \mathcal{L}^h_k z^h) : \mathbb{A}_k : \nabla \hat{\xi}^h \dd x = \int_\Omega \mathcal{L}^h_k z^h \cdot \hat{\xi}^h \dd x
=\int_\Omega \nabla z^h :\Abb_{k}:\nabla \hat{\xi}^h\, \dd x
\end{eqnarray*}
for all $\hat{\xi}^h \in \Xrm^h$; hence, by taking $\hat{\xi}^h = \mathcal{G}^h_k \mathcal{L}^h_k z^h - z^h$ we deduce that $\|\nabla(\mathcal{G}^h_k \mathcal{L}^h_k z^h - z^h )\|^2_{\Lrm^2} = 0$, whereby $\mathcal{G}^h_k \mathcal{L}^h_k z^h = z^h$.

Let us first concentrate on the case when $d=3$. By writing $\xi^h_k:= \mathcal{L}^h_k z^h$, noting~\eqref{eq:proj1},
using the Sobolev embedding $\Wrm^{2,2}(\Omega;\R^m) \hookrightarrow \Wrm^{1,6}(\Omega;\R^m)$, an inverse inequality (recall that $\mathcal{T}^h$ has
been assumed to be quasiuniform), and that by an Aubin--Nitsche duality argument
\[\|\mathcal{G}_k\xi^h_k - \mathcal{G}^h_k\xi^h_k\|_{\Lrm^2}\leq C h^2\|\mathcal{G}_k\xi^h_k\|_{\Wrm^{2,2}} \leq C h^2 \|\xi^h_k\|_{\Lrm^{2}},\]
we have that
\begin{align*}
\|z^h\|_{\Lrm^6}  &= \|\nabla  \mathcal{G}^h_k \mathcal{L}^h_k z^h\|_{\Lrm^6} = \|\nabla  \mathcal{G}^h_k \xi^h_k\|_{\Lrm^6}\\
& \leq \|\nabla  \mathcal{G}_k \xi^h_k\|_{\Lrm^6} + \|\nabla  (I - \mathbb{P}^h)\mathcal{G}_k \xi^h_k\|_{\Lrm^6} +
\|\nabla  (\mathbb{P}^h\mathcal{G}_k - \mathcal{G}^h_k)\xi^h_k\|_{\Lrm^6}\\
& \leq C \|\mathcal{G}_k \xi^h_k\|_{\Wrm^{2,2}} + C \|\mathcal{G}_k \xi^h_k\|_{\Wrm^{1,6}} + C h^{-1 + 3(\frac{1}{6}-\frac{1}{2})}
\|\mathbb{P}^h(\mathcal{G}_k \xi^h_k) - \mathcal{G}^h_k\xi^h_k\|_{\Lrm^2} \\
& \leq C \|\mathcal{G}_k \xi^h_k\|_{\Wrm^{2,2}} + Ch^{-2} (\|\mathbb{P}^h(\mathcal{G}_k \xi^h_k) - \mathcal{G}_k\xi^h_k\|_{\Lrm^2} + \|\mathcal{G}_k\xi^h_k - \mathcal{G}^h_k\xi^h_k\|_{\Lrm^2})\\
& \leq C \|\mathcal{G}_k \xi^h_k\|_{\Wrm^{2,2}} + Ch^{-2} (Ch^2\|\mathcal{G}_k \xi^h_k\|_{\Wrm^{2,2}} + Ch^2\|\xi^h_k\|_{\Lrm^{2}})\\
& \leq C \|\mathcal{G}_k \xi^h_k\|_{\Wrm^{2,2}} + C \|\xi^h_k\|_{\Lrm^{2}}\\
& \leq C \|\xi^h_k\|_{\Lrm^{2}} = C \|\mathcal{L}^h_k z^h\|_{\Lrm^{2}},
\end{align*}
as required.

For $d=2$ the proof is analogous: the $\Lrm^6$ norm is replaced by the $\Lrm^p$ norm for $p \in (1,\infty)$, the embedding
$\Wrm^{2,2}(\Omega;\R^m) \hookrightarrow \Wrm^{1,6}(\Omega;\R^m)$ is replaced by $\Wrm^{2,2}(\Omega;\R^m) \hookrightarrow \Wrm^{1,p}(\Omega;\R^m)$,
$p \in (1,\infty)$, and the exponent $-1 + 3(\frac{1}{6}-\frac{1}{2})$ in the third line of the displayed chain of inequalities above is replaced by $-1 + 2\min\{(\frac{1}{p}-\frac{1}{2}),0\}$.

\subsection*{Acknowledgments}

The authors would like to thank Alexander Mielke, Tom\'a\v{s} Roub\'i\v{c}ek, Ulisse Stefanelli, and Florian Theil for helpful discussions related to this work. Filip Rindler\ gratefully acknowledges the support from the EPSRC Research Fellowship on ``Singularities in Nonlinear PDEs'' (EP/L018934/1); Sebastian Schwarzacher thanks the program PRVOUK P47; Sebastian Schwarzacher and Endre S\"uli\ thank the ERC-CZ program MORE of the Charles University Prague.


\providecommand{\bysame}{\leavevmode\hbox to3em{\hrulefill}\thinspace}
\providecommand{\MR}{\relax\ifhmode\unskip\space\fi MR }
\providecommand{\MRhref}[2]{%
  \href{http://www.ams.org/mathscinet-getitem?mr=#1}{#2}
}
\providecommand{\href}[2]{#2}

\end{document}